\documentclass[11pt]{amsart}

\usepackage[backref,bookmarks=false]{hyperref}

\usepackage{color}

\numberwithin{equation}{section}

\newcommand{\ben}{\begin{enumerate}}
\newcommand{\een}{\end{enumerate}}

\newcommand{\bea}{\begin{eqnarray}}
\newcommand{\ba}{\begin{array}}
\newcommand{\bean}{\begin{eqnarray*}}
\newcommand{\ea}{\end{array}}
\newcommand{\eea}{\end{eqnarray}}
\newcommand{\eean}{\end{eqnarray*}}
\newcommand{\beq}{\begin{equation}}
\newcommand{\eeq}{\end{equation}}
\newcommand{\bthm}{\begin{thm}}
\newcommand{\ethm}{\end{thm}}
\newcommand{\blem}{\begin{lem}}
\newcommand{\elem}{\end{lem}}
\newcommand{\bprop}{\begin{prop}}
\newcommand{\eprop}{\end{prop}}
\newcommand{\bcor}{\begin{cor}}
\newcommand{\ecor}{\end{cor}}
\newcommand{\bdfn}{\begin{dfn}}
\newcommand{\edfn}{\end{dfn}}
\newcommand{\brem}{\begin{rem}}
\newcommand{\erem}{\end{rem}}
\newcommand{\bpf}{\begin{proof}}
\newcommand{\epf}{\end{proof}}
\newcommand{\bfact}{\begin{fact}}
\newcommand{\efact}{\end{fact}}
\newcommand{\bobs}{\begin{obs}}
\newcommand{\eobs}{\end{obs}}
\newcommand{\bexam}{\begin{exam}}
\newcommand{\eexam}{\end{exam}}
\newcommand{\bclaim}{\begin{claim}}
\newcommand{\eclaim}{\end{claim}}

\newtheorem{thm}{Theorem}[section]
\newtheorem{prop}[thm]{Proposition}
\newtheorem{lem}[thm]{Lemma}

\newtheorem{cor}[thm]{Corollary}
\newtheorem{dfn}[thm]{Definition}
\newtheorem{rem}[thm]{Remark}
\newtheorem{fact}[thm]{Fact}
\newtheorem{claim}[thm]{Claim}
\newtheorem{obs}[thm]{Observation}
\newtheorem{exam}[thm]{Example}

\newtheorem{condition}{Condition}
\newtheorem*{condition'}{Condition 2'}

 \newtheoremstyle{claimstyle}%
   {}
   {}
   {\normalfont}
   {}
   {\itshape}
   {.}
   { }
   {\thmnote{#3}}

\theoremstyle{claimstyle}

\alph{enumii} \roman{enumiii}

\def\cA{\mathcal A}             \def\cB{\mathcal B}       \def\cC{\mathcal C}
\def\cH{\mathcal H}             \def\cF{\mathcal F}       
\def\cL{{\mathcal L}}           \def\cM{\mathcal M}        \def\cP{{\mathcal P}}
             \def\cV{\mathcal V}       \def\cJ{\mathcal J}   \newcommand{\J}{\mathcal{J}}
             
                    \def\cO{\mathcal O}
\def\mF{\mathfrak F}

                      \def\R{{\mathbb R}}
\def\C{{\mathbb C}}

\newcommand{\cbar}{\hat{{\mathbb C}} }

\def\a{\alpha}                \def\b{\beta}             \def\d{\delta}
                          
\def\g{\gamma}                \def\Ga{\Gamma}           \def\l{\lambda} 
\def\La{\Lambda}                         
               \def\sg{\sigma}
               \def\th{\theta}           \def\vth{\vartheta}
\def\ka{\kappa}

\newcommand{\ep}{\varepsilon}
\newcommand{\ph}{\varphi}
\newcommand{\al}{\alpha}
\newcommand{\ga}{\gamma}

\def\1{1\!\!1}

\def\and{\text{ and }}

        \def\diam{\text{\rm {diam}}}
\def\dist{\text{{\rm dist}}}

              \def\bu{\bigcup}
\def\({\bigl(}                \def\){\bigr)}

                        \def\^{\tilde}

            \def\sms{\setminus}
\def\sbt{\subset}

\def\sp{\medskip}             \def\fr{\noindent}        
            
\def\ess{{\rm ess}}
\def\D{{\mathbb D}}

\def\${$ \displaystyle }

\def\den{\rho}
\def\shift{\theta}

\newcommand{\pf}{{\mathcal{L}}}
\newcommand{\pfx}{{\mathcal{L}}_x}
\newcommand{\npf}{\mathcal{\hat L}}
\newcommand{\npfx}{{ \mathcal{\hat L}}_x}

\newcommand{\jul}{\mathcal J}

\newcommand{\julx}{\mathcal J_x}
\newcommand{\fatx}{\mathcal F_x}

\newcommand{\sT}{\mathring T}

\begin{document}

\title[Random dynamics of transcendental functions]{Random dynamics of transcendental functions}


\author{Volker Mayer}
\address{Universit\'e de Lille I, UFR de
  Math\'ematiques, UMR 8524 du CNRS, 59655 Villeneuve d'Ascq Cedex,
  France} \email{volker.mayer@math.univ-lille1.fr \newline
  \hspace*{0.42cm} \it Web: \rm math.univ-lille1.fr/$\sim$mayer}

\author{Mariusz Urba\'nski}
\address{Department of
  Mathematics, University of North Texas, Denton, TX 76203-1430, USA}
\email{urbanski@unt.edu \newline \hspace*{0.42cm} \it Web: \rm
  www.math.unt.edu/$\sim$urbanski}

\date{\today} \subjclass{111}

\begin{abstract}

This work concerns random dynamics of hyperbolic entire and meromorphic functions
of finite order and whose derivative satisfies some growth condition at infinity. This class contains
most of the classical families of transcendental functions and goes much beyond. 
Based on uniform versions of Nevanlinna's value distribution theory we first build a thermodynamical formalism
which, in particular, produces unique geometric and fiberwise invariant Gibbs states.
Moreover, spectral gap property for the associated transfer operator along with exponential decay of correlations
and a central limit theorem are shown. 
This part relies on our construction of new positive invariant cones that are adapted to the setting
of unbounded phase spaces. 
This setting rules out the use of Hilbert's metric along with the usual contraction principle. However these
cones allow us to apply a contraction argument 
stemming from Bowen's initial approach. 
 \end{abstract}

\maketitle

Random dynamics is actually a quite active field. An overview can be found in Arnold's book \cite{Arnold98}
and in Kifer and Liu's chapter in \cite{KifPei06}. The first work on random rational functions is due to 
Fornaess and Sibony \cite{FS91}. Related to this is Rugh's paper on random repellers \cite{Rug08}
and Sumi's work on rational semi-groups (see for example \cite{Si97, Si11}).
A complete picture including thermodynamics and spectral gap is contained in \cite{MSUspringer}
which concerns a much wider class of distance expanding random maps, a class originally introduced by Ruelle \cite{Rue78}. Recently random dynamics of countable infinite Markov shifts \cite{DKS08, Stadl10} and graph directed Markov systems \cite{RoyUrb2011} have been treated.
 Here we extend the picture to a situation where the maps are also countable infinite -- to -- one, where the
  phase space is not compact and where in addition there is no Markov structure. 

Given a probability space $(X,\cF , m)$ along with an invertible ergodic transformation $\th: X\to X$,
we consider the dynamics of 
$$f_x^n = f_{\shift ^{n-1} (x)}\circ ... \circ f_x \quad , \quad n\geq 1$$
where $f_x:\C\to \cbar $, $x\in X$, is a family of transcendental functions  depending measurably on $x\in X$.
Like in the deterministic case, the normal family behaviour of $(f_x^n)_n$ splits the plane into two parts
and one is interested in the chaotic part $\J_x$, called fiber Julia set. 
Quite general \emph{transcendental random systems} $f_x:\J_x\to \J_{\shift (x)}$, $x\in X$, are considered in this paper
and, as already has been mentioned among the major difficulties one encounters is that the phase space
$\J_x$ is unbounded and the functions are of infinite degree. The members of such random system, are the fiber maps $f_x:\J_x\to \J_{\shift (x)}$ for each particular $x\in X$.

In the deterministic case, this difficulty has been overcome
in \cite{MyUrb10.1} for a very general class of functions, so called balanced meromorphic functions.
They contain most classical families like all periodic functions (tangent, sine, exponential and elliptic),
functions with polynomial Schwarzian derivative, the cosine-root family and many more
(all these examples are discussed in detail in \cite{MyUrb10.1}).
The key point there was to replace the Euclidean metric by a metric having an appropriate singularity at 
infinity. Once this is done, one can use Nevanlinna's value distribution theory to show that the corresponding
transfer operator is well defined and bounded.
The present paper treats random dynamics generated by the families of functions considered in \cite{MyUrb10.1}.
Again we start with an appropriated choice of metric in order to be able to control the transfer operator.
This time we make use of the uniform versions of Nevanlinna's theorems in Cherry-Ye's book \cite{CherryYe01}
(all needed details of this theory are in the Appendix \ref{sec: nev}).

Then, since we are dealing with random dynamics, measurability of all involved operators, measures and functions
has to be checked.
This point has sometimes been neglected in the literature (see the discussion in \cite{Kif08})
or is the reason for additional assumptions. 
Here we take advantage of Crauel's framework \cite{Cra02} and 
treat measurability very carefully. Moreover, this allows us to have a global, in terms of skew product,
approach which, for example, produces directly measurable families of conditional measures
(see Section \ref{subsec 3.2}). This is in contrast to \cite{MSUspringer} where these objects are constructed fiberwise and then later proven to be measurable.\\[-0.3cm]

Having then good behaving transfer operators and measurability, we can proceed with building the thermodynamical formalism.
As the result, we prove the existence and uniqueness of fiberwise conformal measures and the existence and uniqueness of invariant densities. This gives rise to the existence and uniqueness of fiberwise invariant measures absolutely continuous with respect to the conformal ones (see Theorem \ref{thm Gibbs} and Theorem \ref{thm expo} (1)). \\[-0.3cm]

Contrary to Ionescu Tulcea-Marinescu's theorem
\cite{IoTulMar50} (or its generalization by Hennion \cite{Hennion93}),
the method introduced by Birkhoff \cite{Bir57} and developed further by Liverani \cite{Liv95}, based on positive cones and the Hilbert distance,
can be employed in random dynamics. It especially permits us to obtain the spectral gap property.
But this only does work if the phase spaces are compact. In the present paper this is not the case and so the Hilbert distance is of much less use. Indeed, cones
of functions of finite distance are too small since all of its members must be comparable near infinity.
Fortunately there is a very nice contraction lemma in Bowen's manuscript \cite{Bow75}.
In order to be able to adapt it to the present setting, we first produce, via a delicate construction,
non-standard appropriate invariant cones. Once this is done, the Bowen-like argument is quite elementary. In this sense,
the present work, incidentally, simplifies the deterministic work \cite{MyUrb10.1} which uses 
Marinescu-Ionescu-Tulcea Theorem. \\[-0.3cm]

In conclusion, we get the spectral gap property of Theorem \ref{thm expo} (2). It then almost immediately 
implies the version of exponential decay of correlations in Theorem \ref{thm-correlations}
and the Central Limit Theorem \ref{CLT}.


%
%

\section{Transcendental random systems}

We consider  random dynamics generated by a quite general class of transcendental entire or meromorphic functions of finite order.
As in Arnold \cite{Arnold98}, the randomness is modeled by a measure preserving dynamical system
$(X, \mF , m , \th )$, where $(X, \mF , m)$ is a complete probability space and $\th:X\to X$ an invertible measure preserving ergodic transformation. We do not assume the $\sg$--algebra $\mF$ to be countably generated.
To every $x\in X$ associated is an entire or meromorphic function
$$
f_x:\C\to \cbar .
$$
The order of this function is supposed to be finite and is denoted by $\rho(f_x)$.
For every given $z\in \C$, the map $x\mapsto f_x(z)$ is  assumed to be (at least) measurable
as a map from $(X, \mF )$ to $(\C, \cB)$ where $\cB$ is the Borel $\sg$--algebra of $\C$.
We will often call 
$$
(f_x:\C\to \cbar)_{x\in X}
$$ 
a system or, more fully,  a \emph{transcendental random system} or even a \emph{transcendental random dynamical system}
if it satisfies the following four natural conditions.

\begin{condition}[Common growth of characteristic function]\label{C1}
There are two constants $\rho,  C_\rho >0$ and an increasing function $\omega : [0,\infty) \to [0,\infty)$ satisfying $\lim_{r\to\infty} \log r/\omega (r) =0$ such that
$$
\omega (r) \leq \sT _x (r) \leq C_{\rho} r^\rho \quad \text{for all} \quad
r>0 \text{ and all } x\in X\,.
$$
\end{condition}
\noindent Here, following the standard notation in Nevanlinna theory, we denoted by $\sT_x(r)= \sT(f_x,r)$ the spherical characteristic function of $f_x$. All necessary details on Nevanlinna theory and his fundamental main theorems, in the form most convenient for us,   are collected in \ref{sec: nev}. Appendix. Notice that the right hand side inequality of this condition implies that the orders $\rho (f_x) \leq \rho$
whereas the left hand side is simply a quantitative way of saying that the functions $f_x$ are transcendental ($\sT_f (r)= \cO(\log r)$ means that $f$ is a rational function).

\

In order to study the behavior of the orbits $z\mapsto f_x(z)\mapsto f_{\shift(x)}(f_x(z)) \mapsto ...$ it is natural to use the notation
$$f_x^n:= f_{\shift^{n-1}(x)}\circ ...\circ  f_{\shift (x)}\circ f_x \;\; , \;\; n\geq 1\;.$$
For every $x\in X$ the fiber Fatou set $\fatx$ is the set of all points $z\in \C$ for which there exists a neighborhood $U$ of $z$ on which all the iterates $f_x^n$ are well defined and form a normal family. The complements in the plane, 
$$
\julx:= \C \setminus \fatx,
$$
are called the fiber Julia sets. We would like to mention that the maps
$$
X\ni x\mapsto \fatx
$$
and 
$$
X\ni x\mapsto \julx
$$ 
respectively form open random and closed random sets in the sense of Crauel (see \cite{Cra02}) as defined in Section~\ref{subs_random_obs_and_meas}. This however is not important at the moment. Sometimes, they will be also be denoted by $\J(f_x)$, $x\in X$.
We impose the following normalization which mainly signifies that the Julia set does not accumulate at infinity. Here and throughout the whole paper we will use the notation
$$
\D(z,r):=\{w\in\C:|w-z|<r\} 
$$
and 
$$
\D_T:=\D(0,T)
$$ 
\begin{condition}\label{C2}
There exists $T>0$ such that 
$$\Big(\julx \cap \D_T \Big)\cap f_x^{-1}\left(\jul_{\shift (x)} \cap \D_T\right) \neq \emptyset
\;\; , \;\; x\in X\,.$$
\end{condition}
Let $z_x\in \J_x\cap \D_T \cap f_x^{-1} \left( \J_{\shift (x)}\cap \D_T\right)$. We will see in Lemma \ref{25.2 2}
that these points can be chosen in a measurable way. Consider then the translations 
$T_x(z)= z+z_x$, $x\in X\,.$
They conjugate $(f_x)_{x\in X}$ to  a new system, say $(g_x)_{x\in X}$ which again does depend measurably on $x$
and such that 
$$0\in \J(g_x) \quad \text{and} \quad |g_x(0)|\leq 2T\, , \;\;x\in X.$$ 
Notice that the family of translations $(T_x)_{x\in X}$ and the family of its inverses are equicontinuous since $|z_x|\leq T$, $x\in X$. In \cite{MSU13} families of conjugations with this property are, natuarally, called
bi-equicontinuous and they are important since such conjugations preserve topological features of the dynamics. In particular they preserve corresponding Julia sets whereas general conjugations do not, as can be seen from Example 2.3 in \cite{MSU13}.
In conclusion, up to such a conjugation and by replacing the constant $T$ by $2T$ if necessary, 
we can use the following normalizing requirement instead of Condition \ref{C2}: 
\beq \label{appendix label 1}
0=z_x\in \J_x \quad \text{and} \quad |f_x(0)|\leq T\, , \;\; x\in X.
\eeq

\

A straightforward generalization of the notion of hyperbolicity used
in \cite{MyUrb08, MyUrb10.1} to the random setting is the following.

\bdfn\label{dfn 1}
A transcendental random system $(f_x)_{x\in X}$ is called
\ben
\item \emph{topologically hyperbolic}
if there exists $0<\d_0\leq \frac14$ such that for every $x\in X$, $n\geq 1$ and $w\in \jul_{\shift^n(x)}$
all holomorphic inverse branches of $f_x^n$ are well defined on $\D (w, 2\d_0)$.
\item \emph{expanding} if there exists $c>0$ and $\g >1$ such that 
$$|(f_x^n)'(z)|\geq c\g ^n$$
for every $z\in \julx \setminus f_x^{-n} (\infty)$ and every $x\in X$. a
\item \emph{hyperbolic} if it is both topologically hyperbolic and expanding.
\een
\edfn

As in the papers  \cite{MyUrb08, MyUrb10.1}, dealing with deterministic systems, we will consider hyperbolic systems for which one has some more information about derivatives.

\begin{condition}[Balanced growth condition] \label{C3}
 There are $\al_2 > \max\{0, -\al_1\}$ and $\kappa \geq 1$ such that for every $x\in X$ and 
every $z\in \julx \setminus f_x^{-1}(\infty  )$,
 \begin{equation}\label{eq intro growth}
 \ka^{-1}(1+|z|)^{\al _1}(1+|f_x(z)|)^{\al_2} \leq |f_x'(z)| \leq  \ka(1+|z|)^{\al _1}(1+|f_x(z)|)^{\al_2}
\, . \end{equation}
\end{condition}

\begin{condition} \label{C4}
For every $R>0$ and $N\geq 1$ there exists $C_{R,N}$ such that
$$|\left( f_x^N\right)'(z)|\leq C_{R,N} \quad \text{for all} \quad 
z\in \D_R\cap f_x^{-N} \big( \D_R\big ) \;\; \text{and } x\in X\,.$$

\end{condition}

\brem \label{rema deterministic}
As it is explained in \cite{MyUrb08, MyUrb10.1}, many families naturally 
satisfy the balanced growth condition. For these families and for certain classes of entire functions it turns out that \eqref{eq intro growth} entails their order to be $\a_1+1$.
All other conditions, i.e. Conditions \ref{C1}, \ref{C2} and \ref{C4}, are automatically satisfied in the deterministic case.
\footnote{For some very special examples, the lower bound in Condition \ref{C1} can fail. Notice however that,
if $f$ is not a rational function, then $\sT $ growths faster than $\log r$ and this is exactly the property we really need.}
  Therefore, the present setting is a straightforward generalisation of
the deterministic situation, the only difference being that
$\al_2$ in Condition \ref{C3} is constant  whereas it is allowed to be a bounded function in \cite{MyUrb10.1}.
\erem

Throughout this section and also in the rest of this paper we use some standard 
notations. For example, $a\preceq b$ means that $a\leq c b$ for some constant $c$ which does not depend on the involved variables. We also use $\cV _\d(K)$ for the 
$\d$--neighborhood of $K$ in Hausdorff distance generated by the standard Euclidean metric. 


\

\subsection{Mixing}
We shall need the following mixing property.
\blem \label{lem: mixing}
Let $(f_x)_{x\in X}$ be a hyperbolic transcendental random system. Then,  for all $r>0$ and $R>0$
there exists $N=N(r,R)$ such that 
$$f_x^n(\D(z,r)) \supset \overline \D_R\cap \jul_{\shift^n(x)} \quad 
\text{for every } n\geq N, \;\; z\in \julx\cap\overline \D_R \text{ and } x\in X\,.$$
\elem

\bpf
Suppose to the contrary that there exist $r,R>0$ and arbitrarily large integers
$n\geq N$ such that for some $x_n\in X$ and $z_n\in \J _{x_n} \cap \overline \D_R$ 
there exists a point $$w_n\in \left( \overline \D_R \cap \J_{\shift ^n (x_n)}\right)\setminus f_{x_n}^n (\D(z_n,r)) .$$
Define then $\ph _n :\D \to \cbar$ by $\ph _n (\xi ) =  f_{x_n}^n (z_n+ r\xi ) -w_n$.
Note that the family $(\ph _n)_n$ is not normal at the origin. Consequently, there exist arbitrarily large integers $n$ such that 
$$
\ph _n \left(\D(0,1/2) \right) \cap \D(0, \d )\neq \emptyset\,.
$$
But then, it follows from hyperbolicity and, in particular, from the expanding property that
$$
f_{x_n}^{-n} (\D(w_n , \d))\subset \D(z_n, r)
$$
provided that $n\ge q$ is sufficiently large, where $f_{x_n}^{-n}$ is a 
appropriated holomorphic inverse branch of $f_{x_n}^{n} $ defined on $\D(w_n , \d)$.
 But this contradicts the fact that $w_n\not\in f_{x_n}^{n}(\D(z_n,r))$.
\epf

\section{Transfer operators}
Let $\cC_b(\J_x)$ be the space of continuous bounded real--valued functions on $\J_x$
and $\cC_0(\J_x)$ its subspace consisting of all functions converging to $0$ at $\infty$.
Let $(f_x)_{x\in X}$ be a hyperbolic transcendental random system and define:
$$\pfx g(w) = \sum_{f_x(z)=w}e^{\ph_x (z)}g(z) \; , \quad w\in \jul_{\shift (x)} \; \text{ and } \; g\in \cC_b(\J_x).
$$
This is the associated family of transfer operators with potential $\ph_x:\julx \to \R$.
A natural choice for the potentials is $\ph_x=-t\log|f_x'|$ since usually one can choose the parameter $t$ such that these potentials encode the geometric properties
of the dynamical system. In fact, throughout the paper we do deal only with potentials of this form. However, since $f_x$ is of infinite degree, $\pfx$ is in general not well--defined for such potentials.
One might replace it by its spherical version. Then, at least for $t=2$, $\pfx$ would be well-defined but the new obstacle would then arise, that, except for some special cases, $\pfx$ would not be a bounded operator; to see it the reader is invited just to try and write it down for the exponential family. However, using Nevanlinna theory, we showed in \cite{MyUrb08, MyUrb10.1}, still for the deterministic case, that there is a Riemannian metric, somehow in between the Euclidean and spherical one, conformally equivalent to any of them, such that the transfer operator has all the properties needed for developing the thermodynamic formalism and that this holds for the values of parameter $t$ in a sufficiently large domain, containing in particular the hyperbolic dimension. This method does work as soon as the derivative growth condition, i. e. Condition~3, is satisfied. 

So, suppose that $(f_x)_{x\in X}$ satisfies the Condition \ref{C3}. Then, $\al =\al_1+\al_2 >0$. Given any $t>\frac{\rho}{\al}$
there is $\tau \in (0,  \al_2)$ such that 
\beq\label{eq tau}
t>\frac{\rho}{\hat \tau}>\frac{\rho}{\al}\quad \text{where} \quad \hat \tau = \al_1 +\tau\,.
\eeq
\brem\label{rem tau}
Notice that here $\tau$ can be chosen individually for each $t>\frac{\rho}{\al}$. In particular,
we may suppose that $\a_2-\tau>0$ is arbitrarily small.
\erem
Consider then the Riemannian metric
$$
d\sigma_\tau (z) = \frac{|dz|}{(1+|z|)^\tau}\,.
$$
We denote by $|f_x'|_\tau$ the derivative of $f_x$ with respect to this metric, and, using Condition \ref{C3}, we have,
$$ 
|f_x' (z) |_\tau= |f_x'(z)|\frac{(1+|z|)^\tau}{(1+|f_x(z)|)^\tau}
\asymp (1+|f_x(z)|)^{\al_2-\tau} ( 1+|z|)^{\hat\tau}
$$
for every $z\in \julx \setminus f_x^{-1}(\infty  )$. Denote by $\1$ the function identically equal to $1$ on its appropriate domain. By virtue of the above formula, for all $w\in \jul _{\shift (x)}$ we have, 
\beq \label{2.?}
\pfx \1(w) =\pf _{x,t}\1(w)= \sum_{f_x(z)=w}|f_x' (z) |_\tau^{-t}\leq \frac{\kappa^t  }{(1+|w|)^{(\al_2-\tau ) t}}
\sum_{f_x(z)=w}\big( 1+|z|\big)^{-t\hat\tau}\;.
\eeq

\brem\label{3.3 1} Hyperbolicity of the functions $f_x$ implies that, 
increasing $\kappa$ if necessary, \eqref{2.?} does hold
for all $w$ in the $\d_0$--neighborhoods $ \cV_{\d_0} ( \jul _{\shift (x)} )$ of the Julia sets provided $\d_0>0$ has been chosen sufficiently small.
\erem

\noindent Since the factor with $w$, appearing in the right hand side of \eqref{2.?},
converges to zero as $|w|\to \infty$, and applying also Nevanlinna theory 
(similarly as in \cite{MyUrb10.1}, details can be found in \ref{sec: nev} of Appendix)),  we see that the series in \eqref{2.?} can be uniformly bounded from above. We therefore obtain the following good behavior of these operators $\pfx$.

\bprop \label{prop bounded operators}
For every $t>\frac{\rho}{\hat \tau}>\frac{\rho}{\al}$ there exists $M_0=M_0(t,\tau)>0$ such that for every $x\in X$, we have
\ben
\item \$\sum_{f_x(z)=w}\big( 1+|z|\big)^{-t\hat\tau} \leq M_0$ \ for every \ $w\in \J_{\shift (x)}$,
\item $\|\pfx \|_\infty \leq M_0$ \ and
\item $\pfx\1 (w) \leq M_0 \left( 1+|w|\right)^{-(\al_2-\tau )t} \longrightarrow 0\;$ as $\;|w| \to \infty$.
\een
\eprop

\subsection{Distortion and H\"older functions}  

Koebe's Distortion Theorem (see Theorem 1.3 in \cite{PommerenkeBook}) and elementary calculus give: 

\blem \label{lem koebe}
Given $t,\tau>0$ as in \eqref{eq tau}, there exists $K=K_{t,\tau}>0$ such that,
for every $x\in X$, every integer $n\geq 1$ and every $\psi$, an inverse branch of $f_x^n$ defined on some disk $\D(w,2\d_0)$, $w\in \jul_{\shift^n (x)}$, we have that
$$
\frac{|\psi'(w_1)|_\tau^t}{|\psi'(w_2)|_\tau^t} \leq 1+ K|w_1-w_2|, \quad w_1,w_2\in \D(w,\d_0 )\,.
$$
\elem
\brem
It is reasonable to require that $1/2\leq \tau \leq 2$. This would then imply that the constant $K$ does depend only on the parameter $t$.
\erem

\noindent A Straightforward application of Lemma \ref{lem koebe} gives (remember that $\d_0\leq 1/4$):

\blem \label{lem dist operator} 
There exists $K=K_{t,\tau }$ such that, for every $x\in X$, every integer $n\geq 1$, and every $w\in \jul_{\shift ^n(x)}$, we have that
$$
\frac{\pfx^n\1 (w_1)}{\pfx^n\1 (w_2) }\leq 1+K|w_1-w_2| \quad \text{for all} \quad w_1,w_2\in \D(w,\d_0)\,.
$$
In particular,
$$\pfx^n\1 (w_1)\leq K \pfx^n\1 (w_2) \quad \text{for all} \quad w_1,w_2\in \D(w,\d_0)\,.$$
\elem

\

Let $\cH_\b(\julx)$ be the set of real-valued bounded $\b$--H\"older functions defined on $\julx$.
The $\b$--variation of a function $g\in \cH_\b (\julx )$ is defined to be
\beq\label{l21.9}
v_\b(g):= \sup_{0<|w_1-w_2|\leq \d_0} \left\{\frac{|g(w_1)-g(w_2)|}{|w_1-w_2|^\b}\right\}
\eeq
and 
$$
\|g\|_\b:= v_\b (g) + \|g\|_\infty
$$ 
is the corresponding $\b$--H\"older norm on $\cH_\b (\julx )$. The good distortion behavior established in Lemma \ref{lem koebe} implies the following two-norm inequality
which, in particular, yields invariance of H\"older spaces $\cH_\b(\julx)$.

\bprop \label{prop two norm}
Let $c,\g>0$ be the expanding constants from Definition \ref{dfn 1}.
Then, for every $x\in X$, every integer $n\geq 1$, and every $g\in \cH_\b (\julx)$, we have, with  some $K=K_{t,\tau }>0$ that
$$
v_\b (\pfx ^n g)\leq \|\pfx ^n\|_\infty \Big( \|g\|_\infty +  K(c\g^n)^{-\b} v_\b (g))\Big)\,.
$$
\eprop

\bpf
Let $x\in X$, $n\geq 1$, $g\in \cH_\b (\julx )$ and let $w_1,w_2\in \jul_{\shift^n (x)}$ with $|w_1-w_2|\leq \d_0$. The points $z_1,z_2$ are said to form a pairing if they are respectively preimages of $w_1$ and $w_2$ by the same holomorphic inverse branch of $f_x^n$. With this convention,
$$
|\pfx g (w_1)-\pfx g (w_2)|\leq  I+II,
$$
where
$$
II=:\sum_{z_1,z_2\;  pairing} |(f_x^n)'(z_2)|_\tau ^{-t} |g(z_1)-g(z_2)|
\leq \|\pfx^n\|_\infty v_\b (g) (c\g ^n)^{-\b} K^\b |w_1-w_2|^\b$$
and 
\begin{equation*}
    \begin{split} 
I=& \sum_{z_1,z_2\;  pairing} \Big||(f_x^n)'(z_1)|_\tau ^{-t}- |(f_x^n)'(z_2)|_\tau ^{-t} \Big|g(z_1)\\
\leq &\|g\|_\infty  \sum_{z_1,z_2\;  pairing} |(f_x^n)'(z_1)|_\tau ^{-t} \left|1- \frac{|(f_x^n)'(z_2)|_\tau ^{-t}}{|(f_x^n)'(z_1)|_\tau ^{-t}} \right| \\
\leq &\|g\|_\infty  \| \pfx^n\|_\infty K |w_1-w_2|,
    \end{split} 
\end{equation*}
where the last inequality results from Lemma~\ref{lem koebe}.
It suffices now to combine the above estimates of both terms $I$ and $II$. 
\epf


%
%


\

\section{Random Gibbs states}

In this section we establish the following key result which the rest of the paper relies on.

\bthm \label{thm Gibbs}
Let $(f_x)_{x\in X}$ be a hyperbolic transcendental random dynamical system satisfying Conditions \ref{C1}-\ref{C4}. Fix $t>\rho / \al$. Then there exists a random Gibbs measure $\nu$ with disintegrations $\nu_x\in \cP\cM ( \J_x)$, $x\in X$,
and a measurable function $\l : X\to (0, \infty ) $ such that 
\beq \label{eq Gibbs}
\cL^* _x \nu_{\shift (x)}=\l_x \nu _x \quad  \text{for } m \text{ a.e. } x\in X\,. 
\eeq
Moreover, there exists a constant $C \geq 1$ such that \ $C^{-1}\leq \l_x \leq C$ for $m$--a.e. $x\in X$.
\ethm

\noindent Yes, we have not defined random measures yet. Roughly speaking, this concept means that the family of probability measures $(\nu_x)_{x\in X}$ is measurable and to integrate one integrates first aginst the fibers over points $x\in X$ and then against the measure $m$ on $X$. In order to get measurability of $\l$ and $\nu$, unlike to the previous sections, it is much better now to consider the global skew product map
$$
(x,z)\mapsto F(x,z):=(\shift (x) , f_x (z)),
$$
the global transfer operator, and the associated global Julia set 
\beq \label{eq J}
\J =\bigcup_{x\in X} \{ x\} \times \J_x \subset X\times \C,
\eeq
along with the measurable structure of $\J$ induced by the $\sg$--algebra $\mF\otimes \cB$ of $X\times \C$
where $\cB$ is the Borel $\sg$--algebra of $\C$. The advantage is that then one can use the framework of random sets and measures described by Crauel in \cite{Cra02}. We now first present this framework along with some applications, and then prove Theorem \ref{thm Gibbs}.

\subsection{Random observables and measures}\label{subs_random_obs_and_meas}
Let us recall first that $(X,\mF,m)$ is an arbitrary complete probability space with the spaces $X$, the $\sg$--algebra on $X$, and a complete probability measure $m$ on $(X,\mF)$. Following Definition~2.1 in \cite{Cra02} we say that a function
\beq\label{1mu2015_01_30}
X\ni x\mapsto C_x,
\eeq
ascribing to each point $x\in X$ a closed subset $C_x$ of $\C$, is called a random closed set if for each $z\in \C$ the function
$$
X\ni x\mapsto \dist(z,C_x)\in\R
$$
is measurable. Since the probability measure $m$ on $X$ is assumed to be complete, being a closed random set precisely means (see Proposition~2.4 in \cite{Cra02}) that the union
\beq\label{2mu2015_01_30}
C:=\bu_{x\in X}\{x\}\times C_x \  \text{ (all sets $C_x$ are assumed to be closed)}
\eeq
is a measurable subset of $X\times\C$. This natural bijection entitles us to speak and to refere to measurable subsets of $X\times\C$ with closed sections along $\C$ also as closed random sets. 
We will therefore frequently refer to both the functions as in \eqref{1mu2015_01_30} and the sets of \eqref{2mu2015_01_30} as random sets.

\sp A random closed set $X\ni x\mapsto C_x$ is called a random compact set if all sets $C_x$, $x\in X$, are compact (in $\C$). A function $X\ni x\mapsto V_x$ is called a random open set if the function $X\ni x\mapsto \C\sms V_x$ is a
closed random set. 

\sp A particularly important feature of closed random sets is that they allow us to use a Measurable Selection Theorem, namely Theorem 2.6 in \cite{Cra02}.
This theorem asserts that for any closed random set $C\subset \mF\otimes \cB$
there exists a countable family of measurable functions $(c_n:X\to\C)_{n\geq 0}$ such that for $m$--a. e. $x\in X$, 
\beq \label{25.2 3} 
C_x = \overline{\left\{c_n(x):n\geq 0\right\}}. 
\eeq
We shall prove the following.

\blem\label{25.2 1}
The global Julia set $\J$ is a closed random set.
\elem

\bpf For every $x\in X$ denote by $\cC_{f_x}$ the set of all critical points of $f_x$, i.e.
$$
\cC_{f_x}=\{z\in\C:f_x'(z)=0\}.
$$
Set 
$$
\cO_{n,x}^+ = f_x^n(\cC_{f_x})\cup f_{\shift (x)}^{n-1}(\cC_{f_{\shift (x)}})\cup...\cup
f_{\shift ^{n-1}(x)}(\cC_{f_{\shift^{n-1} (x)}}).
$$ 
and then
$$
\cP_x =\overline{\bigcup_{n\geq 1} \cO_{n,\shift^{-n}(x)}^+}\cap \C\,.
$$
Since each set $\cC_{f_x}$ is countable and its elements vary measurably with $x\in X$, Proposition~2.9 in \cite{Cra02} assures us that $\cP =\bigcup_{x\in X} \{x\}\times \cP_x$ is a closed random set. Hence, still by Proposition 2.4 in \cite{Cra02}, 
$\cP^{\d} =\bigcup_{x\in X}\{x\}\times \overline{\cV_{\d}(\cP_x)}$
is also a closed random set, where $\d=\d_0/2$. It now follows from Proposition~2.9 in \cite{Cra02} that $C= \bigcup_{x\in X} \{x\}\times C_x$, with $C_x= \C\setminus \cV_{\d}(\cP_x)$, is a closed random set. The Measurable Selection Theorem, Theorem~2.6 in \cite{Cra02}, thus applies, and, as $\cV_{\d} (\J_x ) \subset C_x$ for all $x\in X$ by hyperbolicity of $(f_x)_x$, this theorem yields measurable maps $(c_k:X\to \C)_{k\ge 0}$ such that $c_k(x)\in \cV_{\d} (\J_x )$ for every $k\ge 0$, and moreover, for $m$--a.e. $x\in X$, $\overline{\{c_k(x):k\ge 0\}}\supset\cV_{\d}(\J_x )$. By definition of $C$, all holomorphic inverse branches of $f_x^n$
are well--defined in the $\d_0$ neighborhoods of all points $c_{\shift^n (x)}$.
Fix $1< \eta < \g$ where $\g $ is the expanding constant coming from Definition~\ref{dfn 1}. We call a holomorphic inverse branch
$f_{x,*}^{-n}$ of $f_x^n$, defined on $\D(c_{\shift ^n (x)},\d )$, \emph{shrinking}, if 
$|(f_{x,*}^{-n})'(c_{\shift ^n (x)})|\leq \eta ^{-n}$. It is now easy to check that
$$
\J_x= \bigcap _N \overline{\bigcup_{n\geq N} \bigcup_{* \; \text{shrinking}}f_{x,*}^{-n} (\{c_k(\shift^n(x):k\ge 0\})}.
$$
This shows that $\J$ is a closed random set.
\epf

\blem\label{25.2 2}
If Condition \ref{C2} holds, then there is a measurable choice
$$
x\mapsto z_x\in \julx \cap \overline\D_T \cap f_x^{-1}\left(\jul_{\shift (x)} \cap \overline\D_T\right)\,.
$$
\elem 

\bpf
Since, by Lemma \ref{25.2 1}, $\J$ is a closed random set, the sets with fibers $\J_x \cap \overline \D_T$
and $f_x^{-1}(\J_{\shift(x)}\cap \overline \D_T)$ are both closed random sets. The intersection of these closed random sets is again a closed random set, and Condition  \ref{C2} implies that every fiber of this intersection is not empty. 
Therefore, again by the Measurable Selection Theorem, there exist a measurable map $z$ such that $z_x\in \julx \cap \overline\D_T \cap f_x^{-1}\left(\jul_{\shift (x)} \cap \overline\D_T\right)$ for a.e. $x\in X$.
\epf

We now introduce random observables. We recall from \cite{Cra02} that a function $g:\J \to \R$, \ $(x,z)\mapsto g_x(z)$, is called random continuous if $g_x\in \cC_b(\J_x)$ for all $x\in X$, the function $x \mapsto \| g_x \|_\infty$ is measurable and, moreover, 
$m$--integrable. The vector space of all such functions is denoted by $\cC_b (\J )$. It becomes a Banach space when equipped with the norm
$$
|g| = \int _X  \| g_x \|_\infty\, d m(x).
$$
We need more special functions.

\bdfn\label{25.2 5}
A random continuous function $g:\J \to \R$, $(x,z)\mapsto g_x(z)$, is said to vanish at infinity if 
$$
\lim_{z\to\infty}g_x(z)=0
$$
for $m$-a.e. $x\in X$.  The vector space of all such functions is denoted by $\cC_0(\J )$. It is a closed subspace of $\cC_b (\J )$ and inherits the norm $|\cdot|$ from $\cC _b(\J )$. Thus, it becomes a Banach space on its own. 
\edfn

\bdfn\label{25.2 5h}
Random $\b$--H\"older observables are defined to be all the functions $g\in \cC_b(\J )$ such that  $g_x\in \cH_\b (\J_x )$ and such that $x\mapsto \| g_x \|_\b$ is integrable.
This space is denoted by $\cH_\b (\J)$ and equipped with the norm  
$$
| g|_\b = \int _X  \| g_x \|_\b \, d m(x) \, .
$$
\edfn

Consider now the global transfer operator
$\cL $ defined by $$(\cL g)_{x}(w) = \cL_{\shift^{-1}(x)}g_{\shift^{-1}(x)} (w)\; ,\quad (x,w)\in \J\, .$$

\blem \label{25.2 4}
If $g\in \cC_b (\J )$, then $\cL g$ is measurable.
\elem

\bpf
First of all, it suffices to establish measurability of $\cL g$ restricted to measurable sets of the form
$$
E_w =\J \cap (X\times \D(w, \d/2 ))\; , \quad w\in \J_x \; , \;\; x\in X\,.
$$
So, let $(x,w)\in \J$. Since $\J$ is a closed random set, it follows from Proposition~2.4 in \cite{Cra02}) that the 
set
$$
Y=\{y\in X:\J_y \cap \D( w, \d /2 ) \neq \emptyset\}
$$ 
is measurable. Notice that, by definition of $Y$ and by hyperbolicity of $(f_x)_x$, the function $\cL g$ is in fact well defined on $Y\times  \D( w, \d /2 )$. Therefore, we can consider the map $h:X\times \D(w, \d/2 ) \to \R$ defined by
$$
h_y (z)
=\begin{cases}
(\cL g)_y (z) \  &\text{ if }\; y\in Y       \\
            0 \  &\text{ if }\; y\not\in Y.
\end{cases}
$$
Obviously, to show that $\cL g _{| E_w}$ is measurable, it suffices to establish measurability of $h$.
Also, since $Y$ is measurable and $h\equiv 0$ on $Y^c \times \D(w, \d/2 )$ it suffices
to show that $h$ restricted to $Y \times \D(w, \d/2 )$ is measurable and, by virtue of Lemma 1.1 in \cite{Cra02}, in order to prove this, it suffices to show that for every $y\in Y$ the map $\D(w, \d/2 )\ni z\mapsto h_y(z)$ is continuous and that the map $Y\ni y\mapsto h_y(z)$ is measurable for every $z\in \D(w, \d/2 )$.
The continuity for fixed $y\in Y$ is obvious. So we are left to show measurability of
$$
Y\ni y\mapsto h_y (z) = (\cL g )_y (z) = \cL_{\shift^{-1}(y)}g_{\shift^{-1}(y)} (z)
$$
for every fixed $z\in \D(w, \d/2 )$. 
Let $z\in \D(w, \d/2 )$. Then the set $C= \bigcup_{y\in Y} \{y\} \times f_{\shift^{-1}(y)} ^{-1} (z)$
is a closed random set with discrete fibers. Therefore, the Selection Theorem yields the 
existence of countably many measurable functions $(c_n:)_{n=1}^\infty$ such that 
$\{c_n(y):n\ge 1\}= f_{\shift^{-1}(y)} ^{-1} (z)$ for $m$--a.e. $y\in Y$.
Consequently,
$$
h_y (z)
= \cL_{\shift^{-1}(y)}g_{\shift^{-1}(y)} (z) 
= \sum _{n\ge 1} |f'_{\shift^{-1}(y)}(c_n(y))|_\tau ^{-t} g_{\shift^{-1}(y)} (c_n(y)) \quad
\text{for } m-\text{ a.e. }  y\in Y.
$$
This proves the desired measurability.
\epf

\noindent Combining Lemma \ref{25.2 4} with Proposition \ref{prop bounded operators} leads to the following.
\bprop\label{25.2 6}
The transfer operator $\cL$ acts continuously on both $\cC_b (\J )$ and $\cC_0(\J )$.
\eprop

Let $\cM_m(\J)$ be the space of all real-valued signed finite Borel measures $\nu$ on $\J$ such that $\nu\circ\pi_X^{-1}$ is absolutely continuous with respect to $m$ and,
if $(\nu_x)_{x\in X}$ is the corresponding disintegration of $\nu$, then the map 
$$
X\ni x\mapsto ||\nu_x||
$$ 
is in $L^\infty (X)$, i.e. 
$$
||\nu||_\infty:=\ess\!\sup_x ||\nu_x|| <+\infty,
$$
where $|\nu_x|$ is the variation of $\nu_x$ and $||\nu_x||:=|\nu_x|(\1)$ is the total variation of $\nu_x$. 

\sp On the other hand, any function
\beq\label{1mu_2015_01_31}
X\ni x\mapsto \nu_x
\eeq
taking respective values in the spaces of signed measures on $\cJ_x$ such that the function
$$
X\ni x\mapsto \int _{\J_x} g_x\, d\nu_x
$$
is measurable for all functions $g\in\cC_b (\J)$, and $||\nu||_\infty<+\infty$, gives rise to an element $\nu\in\cM_m(\J)$ via the following integration formula. For every  $g\in \cC_b (\J )$:
\beq\label{1mu_2015_02_03}
\nu (g) =\int _X \int _{\J_x} g_x\, d\nu_x dm(x)
\eeq
This enables to speak of elements of $\cM_m(\J)$ either as of appropriate measures or, equivalently, of functions from \eqref{1mu_2015_01_31} as described above.

\sp We denote by $\cM_m^+(\J)$ the subset of $\cM_m(\J)$ consisting of non-negative measures. We further want to single out one particular subspace of $\cM_m(\J )$. This subspace will be essential in the sequel. In its definition, stated below, $\pi_X :X\times\C \to X$ and $\pi_\C :X\times\C \to \C$ are the usual respective projections
$$
\pi_X(x,z)=x \  \  \text{{\rm and }} \  \pi_\C(x,z)=z.
$$
\bdfn\label{random measures}
A measure $\nu \in \cM_m^+(\J )$ with marginal $m$, i.e. such that $\nu\circ\pi_X^{-1}= m$, is called a random measure. In other words this means that $\nu\in \cM_m(\J )$ and the corresponding  disintegrations $(\nu_x)_{x\in X}$ belong to respective spaces $\cP(\J_x )$ of Borel probability measures on $\cJ_x$ for all $x\in X$. The subspace of $\cM_m (\J )$ consisting of all random measures will be denoted by $\cP_m( \J)$. In accordance with the discusion above we identify random measures with functions 
\beq\label{2mu_2015_01_31}
X\ni x\mapsto \nu_x\in\cP(\J_x)
\eeq
such that the function
$$
X\ni x\mapsto \int _{\J_x} g_x\, d\nu_x
$$
is measurable for all functions $g\in\cC_b (\J)$. We also refer to such functions (of \eqref{2mu_2015_01_31}) as random measures.
\edfn

Random measures, as defined in Crauel's book \cite{Cra02}, are measures on the set $X\times \C$. But here we are only interested
in the subclasses  $\cP_m(\J )$ and $\cM_m(\J )$ and they are respective measures in $\cM_m:=\cM_m (X\times \C)$ and $\cP_m:=\cP_m(X\times \C)$ with support in $\J$.

\sp The key concept pertaining to random measures is that of narrow topology, which is a version of weak convergence. Namely, if $\La$ is a directed set, then a net $\(\nu^\a\)_{\a\in\La}$ in $\cM_m$ is said to converge to a random measure $\nu\in \cM_m$ if
$$
\lim_{\a\in\La}\nu^{\a} (g) =\nu (g) \quad \text{for every $g\in \cC_b(X\times\C)$. }
$$ 
This concept of convergence defines a topology on $\cM_m$ called in \cite{Cra02} the narrow topology. The narrow topology on $\cM_m(\cJ)$ is the one inhereted from the narrow topology on $\cM_m$. Since $\1_A$, the characteristic function of $A$, any closed random set in the complement of  $\cJ$, belongs to $\cC_b(X\times\C)$, we have that
$$
\nu(\1_A)=\lim_{\a\in\La}\nu^{\a} (\1_A) =0
$$ 
for any net $\(\nu^\a\)_{\a\in\La}$ in $\cM_m(\cJ)$ converging to a random measure $\nu\in \cM_m$. This means that then $\nu\in \cM_m(\cJ)$, leading to the following.

\bprop\label{p2_2015_01_17}
$\cM_m(\cJ)$ is a closed subset of $\cM_m$ with respect to the narrow topology on $\cM_m$.
\eprop

\bobs
Of course $\cP_m(\cJ)$ is a closed subset of $\cM_m$.
\eobs

\fr Recall that a subset $\Ga$ of $\cM_m$ is bounded if 
$$
\sup\{||\nu||_\infty:\nu\in\Ga\}<+\infty.
$$
Recall from \cite{Cra02} that a subset $\Ga$ of $\cM_m$ is called tight if its projection $\Ga\circ\pi_{\C}^{-1}$ on $\C$ is a tight subset of Borel probability measures on  $\C$, the latter (commonly) meaning that for every $\ep>0$ there exists a compact set $K_\ep\sbt\C$ such that $\nu\circ\pi_{\C}^{-1}(K_\ep^c)\le \ep$ for all $\nu\in\Ga$. For us, the crucial property of narrow topology is that of Prohorov's Compactness Theorem (Theorem 4.4 in \cite{Cra02}) which asserts that a bounded subset $\cM \subset \cM_m$ is relatively compact with respect to the narrow topology if and only if it is closed and tight. Along with Proposition~\ref{p2_2015_01_17}, this entails the following.

\bthm[Prohorov Compactness Theorem of Crauel]\label{theo crauel prohorov}\label{t3_2015_01_17}
A bounded subset $\Ga\subset \cM_m(\cJ)$, in particular, any subset $\Ga\subset \cP_m(\cJ)$, is relatively compact with respect to the narrow topology if and only if it is tight. Furthermore, it is compact if and only if it is tight and closed.
\ethm

\fr The most relevant theorems about tightness  are these (see Proposition~4.3 in \cite{Cra02}):

\bprop\label{p4_2015_01_17}
A subset $\Ga\subset \cP_m(\cJ)$ is tight if and only if for every $\ep>0$ there exists a random compact set $X\ni x\mapsto K_x$ such that $K_x\sbt\cJ_x$ for all $x\in X$ and
$$
\int_X\nu_x(K_x)\, dm(x)\ge 1-\ep
$$
for all $\nu\in\Ga$.
\eprop
\fr As an immediate consequence of this proposition, we get the following. 

\bcor\label{c5_2015_01_17}
Let $\Ga$ be a susbset of $\cP_m(\cJ)$. Suppose that for every $\ep>0$ there exists a random compact set $X\ni x\mapsto K_x$ such that $K_x\sbt\cJ_x$ for all $x\in X$ and
$$
\nu_x(K_x))\ge 1-\ep
$$
for all $\nu\in\Ga$ and for $m$-a.e. $x\in X$. Then the set $\Ga$ is compact. 
\ecor 


\subsection{Random Gibbs states} \label{subsec 3.2}
The aim of this section is to prove Theorem~\ref{thm Gibbs} and to provide some useful estimates. Each measure $\nu_x$, $x\in X$, from Theorem~\ref{thm Gibbs} will be called  a fiber \emph{Gibbs state} or fiber \emph{conformal measure}. The function
$$
X\ni x\mapsto \nu_x
$$
will be, accordingly, called a random Gibbs state or a random conformal measure. Likewise, in accrodance with the above, the global measure (see Definition~\ref{random measures} and formula \eqref{1mu_2015_02_03}), will be called a random Gibbs state or a random conformal measure.


\sp From the invariance relation \eqref{eq Gibbs} it follows that $\l_x = \int \cL_x \1 \, d\nu_{\shift (x)}$ and so we look for measures
$(\nu_x)_{x\in X}$ that are invariant under the map $\Phi : \cP_m(\J) \to \cP_m(\J )$ whose fiber maps
$ \Phi _x : \cP(\jul_{\shift (x)} ) \to \cP(\jul_{x} )$ are defined by
\beq 
\label{3.2 1}  \Phi_x (\nu_{\shift (x)})= \frac{\cL_x^*  \nu_{\shift (x)}}{\cL_x^*  \nu_{\shift (x)} (\1)}
= \frac{\cL_x^*  \nu_{\shift (x)}}{  \nu_{\shift (x)} (\cL_x \1)}\;. 
\eeq 
We want to obtain these measures in the usual way by employing Schauder--Tychonoff's Fixed Point Theorem.
But, since the sets $J_x$ are unbounded, this can be done only if a convex compact and $\Phi$--invariant space of probability measures is found, and if in addition $\Phi$ acts continuously on this space. Towards this end, consider
\beq \label{3.2 2}
\cM=\cM (R_0,  \ep ) := \Big\{ \nu = (\nu_x)_{x\in X}\in \cP_m(\J ): \text{(a) and (b) hold} \Big\},\quad \text{where }
\eeq
\begin{itemize}
\item[(a)]  $\nu_x(\overline {\D} _{R_0} ) \geq \frac12$ and
\item[(b)] $\nu_x(\overline {\D} _{R}^c ) \leq \frac{1}{R^\ep}$ for every $R\geq R_0$
\end{itemize}

\vskip 2.00mm
\noindent are required to hold for $m$--a.e. $x\in X$.

\vskip 2.00mm \noindent
For any $x\in X$, define further 
$$
\cM_x = \{\nu_x: \text{ (a) and (b) hold}\}.
$$
Clearly, property (b) implies (a) with some $R_0$ sufficiently large. The significance of (a) is to specify some radius $R_0$. Invoking Corollary~\ref{c5_2015_01_17} and Theorem~\ref{theo crauel prohorov}, we obviously, we have the following.

\blem \label{3.2 3}
The set $\cM$ is convex, closed and tight, hence compact.
\elem

\noindent We shall prove the following.

\bprop \label{3.2 4}
There are $R_0, \ep>0$ such that $\cM=\cM (R_0,  \ep )$ is invariant under the map  $\Phi=(\Phi _x)_{x\in X}$
defined in \eqref{3.2 1}, i.e. $\Phi_x (\cM_{\shift (x)}) \subset \cM_x$ for all $x\in X$.
\eprop

\noindent In order to establish this result we first need two lemmas.

\blem \label{3.2 5}
For every $0<a< \hat\tau t -\rho $ there exists $M_a$ such that for all $x\in X$ and all $R\ge 1$,
$$
\cL_x\1_{\overline \D _R^c} (w) 
\leq \frac{M_a}{R^a} \quad , \quad 
w\in \J_{\shift (x)}.
$$ 
\elem

\bpf
Given $a\in (0, \hat\tau t -\rho)$, let $b= b(a)>0$ such that $\hat \tau t = a + \rho +b$.
Then,
\begin{align*}\cL_x\1_{\overline \D _R^c} (w) 
&\leq \kappa^t (1+|w| ) ^{-(\a _2 -\tau)t}
 \sum_{z\in f_x^{-1}(w)\cap {\overline \D _R^c}}(1+|z|)^{-t\hat \tau}\\
&\leq \kappa^t  \sum_{z\in f_x^{-1}(w)\cap {\overline \D _R^c}}|z| ^{-a}|z|^{-(\rho +b)}
\leq \frac{M_a}{R^a},
\end{align*}
where, the first inequality fpollows from \eqref{2.?}, while the last one, with some constant $M_a<\infty$, is a consequence of Proposition~\ref{prop bounded operators} with $\hat \tau t$ replaced by $\rho+b$.
\epf

\blem \label{3.2 6}
There exists $\tilde R_0>0$ and, for every $R\geq  \tilde R_0$,
 such that, for some constant $c>0$,
 $$\cL_x\1 (w) \geq c R^{-(\al_2-\tau) t} 
 8\log R \, r_R^{-\hat \tau t } \quad \text{for every} \quad 
   w\in \J_{\shift (x)}\, ,\;\; |w|\leq R\, ,x\in X
$$
where $r_R= \omega^{-1} \big( 8 \log R \big)$.
\elem

\bpf
This proof relies heavily on Nevanlinna's theory and especially on Theorem \ref{usmt}.
The notation used in it is explained in  \ref{sec: nev}' Appendix.
But first we need some preliminary observations.

Since Condition \ref{C2} holds we may assume that $0\in \J_x$ and $|f_x(0)|\leq T$ for  all $x\in X$.
Let $\Omega_x$ be the connected component of $f_x^{-1} \Big( \D (f_x(0) , \d_0 )\Big)$ that contains $0$.
Since $0\in \J_x$ we can use the expanding property of Definition \ref{dfn 1} along with Condition \ref{C4}
in order to get
\beq \label{4.2 1}
c\ga \leq |f_x' (0)|\leq C_T \quad \text{for every} \quad x\in X\,.
\eeq
Koebe's $1/4$--Theorem (see \cite{Hille} for the most common modern source of it proof) applies and, together with \eqref{4.2 1}, implies that 
$$\Omega_x\supset \D\left(0, \frac14 |f_x'(0)|^{-1}\d_0 \right)\supset \D\left(0,s\right)
\quad \text{where} \quad s= \frac{\d_0 }{4C_T} \,.$$
Let from now on $x\in X$ and $w\in \cJ_{\shift (x)}$ with $|w|\leq R$.

\

\noindent
{\emph{Case 1:}} Suppose that $f_x^{-1} (\D(w,\d_0 )) \cap  \D _s \neq \emptyset$, i.e.
that there exists $z'\in \D(0,s)$ with $w'=f_x(z')\in \D(w,\d_0 )$.
Then
$$\cL_x \1 (w)\geq \frac1K \cL_x \1 (w')   
 \geq  \frac1K |f_x' (z')|^{-t}_\tau =
 \frac{|f_x' (z')|^{-t}}{K}\frac{(1+|z'|)^{-\tau t}}{(1+|w'|)^{-\tau t}}
 \geq \frac{(KC_T)^{-t}}{K}(1+s)^{-\tau t}$$
by Lemma \ref{lem dist operator}, Lemma~\ref{lem koebe}, and \eqref{4.2 1}. Hence, in this case there is a uniform lower bound
for $\cL_x\1 (w)$.

\

\noindent
{\emph{Case 2:}} Suppose that $f_x^{-1} (\D(w,\d_0 )) \cap  \D _s = \emptyset$. Then we have to use the 
uniform SMT (Theorem \ref{usmt}) and, in order to do so, first to verify its assumptions. It follows from \eqref{4.2 1} that
$$
\frac{c\ga }{1+T^2}\leq f^\#_x(0)=\frac{|f'_x(0)}{1+|f_x(0)|^2}\leq C_T\quad , \quad x\in X\,.
$$
In other words, the assumption (1) of Theorem \ref{usmt} holds with 
$L=\max \left\{ C_T, \; \frac{1+T^2}{c\ga }\right\}.$
Assumption (3) is exactly the uniform growth condition of the characteristic functions in Condition \ref{C1}.
It remains to choose appropriate points $a_j$. Let 
$a_1,a_2,a_3\in \D( w, \d_0)$ be any points such that
$ |a_i-a_j|\geq \frac{ \d_0}{3}$ for $i\neq j$. Notice that  
$$ 
f_x(0)\not\in\{a_1,a_2,a_3\} \quad \text{and} \quad f_x^{-1} (a_j) \cap \D_s = \emptyset\,.
$$
We need the following simple estimate:
\begin{align*}
\mathring D(a_1,a_2,a_3)=-\log \prod_{i\neq j} [a_i,a_j] +2\log 2\leq &
\log \left(1+|w|^2 \right)+ \log \frac{12}{\d_0} +2\log 2\\
\leq& 2\log R + \log \frac{12}{\d_0} +3\log 2\,.
\end{align*}
Theorem \ref{usmt} of Appendix gives now the following inequality:
\begin{align*}
\sum_{j=1}^3 N(a_j ,r) 
&\geq  \sT _x (r) -b_6 -6\rho \log r -2 \log R -
\log \frac{12}{\d_0} -3\log 2\\
&=  \sT _x (r) -\tilde b_6 - 6\rho \log r -2 \log R \\
&=  \sT _x (r) \left( 1- \frac{\tilde b_6 +6\rho \log r +2 \log R}{\sT _x (r) } \right) .
\end{align*}
Remember that $\sT_x(r) \geq \omega (r) $, that $\lim_{r\to\infty} \log (r) / \omega(r)=0$ (Condition \ref{C1}) and that $R\geq \tilde R_0$. If we define $r_R:=\omega^{-1} \big( 8 \log R\big)$, then
$$\frac{2\log R}{\omega (r)}\leq \frac14 \quad \text{for every} \quad r\geq r_R\,.$$
Therefore, for every $R\geq \tilde R_0$ and provided that $\tilde R_0$ is sufficiently large, we have:
$$
\frac{\tilde b_6 +6\rho \log r +2 \log R}{\omega (r) } \leq \frac12 \quad \text{for every} \quad r\geq r_R.
$$
This implies that
\beq \label{4.2 3}
\sum_{j=1}^3 N(a_j ,r) \geq \frac12  \sT _x (r)\geq \frac12 \omega (r)
\geq \frac12 \omega (r_R) \quad \text{for every} \quad  r\geq r_R \,.
\eeq
We can now conclude the proof of our lemma. Indeed, Lemma~\ref{lem dist operator},  the lower bound in Condition~\ref{C3} 
and the fact that $f_x^{-1} (a_j) \cap \D_{s} =\emptyset$ imply,
for every $j=1,2,3$,
$$
\begin{aligned}
K\cL_x\1 (w) 
&\geq \cL_x \1 (a_j) 
\succeq (1+|a_j|)^{-(\al_2-\tau) t} \sum_{f_x(z)=a_j} (1+|z|)^{-\hat \tau t} \\
&\succeq R^{-(\al_2-\tau) t}  \sum_{f_x(z)=a_j} |z|^{-\hat \tau t}\,.
\end{aligned}
$$
A standard argument  (see \ref{sec: nev}. Appendix here, \cite{Nevbook74} or Chapter 3 of \cite{MyUrb10.1}) and \eqref{4.2 3} shows that
\begin{align*}
 \sum_{j=1}^3\sum_{f_x(z)=a_j} |z|^{-\hat \tau t} = (\hat \tau t)^2 \int _s^\infty \frac{ \sum_{j=1}^3N(a_j,r)}{r^{\hat \tau t +1}}dr
\succeq \omega (r_R) \int _{r_R}^\infty \frac{dr}{r_R^{\hat \tau t+1 }}
\succeq \omega (r_R) r_R^{-\hat \tau t }
\,.
\end{align*}
Finally, there exists $\tilde R_0>0$ and $c>0$ such that
$$3K\cL_x\1 (w) \geq \sum_{j=1}^3  \cL_x \1 (a_j)  \geq c R^{-(\al_2-\tau) t} 
 \omega (r_R) r_R^{-\hat \tau t }$$
for every $R\geq \tilde R_0$ and $w\in \cJ_{\shift (x)}$, $|w|\leq R$, $x\in X$.
\epf

\

\begin{proof}[Proof of  Proposition \ref{3.2 4}.]
Let $\nu _{\shift (x)}\in \cM_{\shift (x)}$. We have to show that there are constants $R_0,\ep$
that do not depend on $x\in X$ such that $\Phi_x(\nu _{\shift (x)})\in \cM_x$, i.e. that the properties (a), (b)
of \eqref{3.2 2} are satisfied. Let $\tilde R_0$ be the number given in Lemma \ref{3.2 6},
suppose that $R_0\geq \tilde R_0$
and let $R\geq R_0$.

We have to choose the constant $a\in (0, \hat \tau t -\rho)$ from Lemma \ref{3.2 5} and $\tau \in (0, \a_2)$.
Let $a=\frac12 (\hat \al t -\rho )$ and, according to Remark \ref{rem tau}, we may choose $\tau$
sufficiently close to $\a_2$ such that $a< \hat \tau t -\rho$ and
$b= \frac{a}{2} - (\a_2-\tau ) t >0$. Lemma \ref{3.2 5} implies  
$$\cL_x^* \nu  _{\shift (x)} (\overline \D_R^c) = \int \cL_x \1 _{\overline \D_R^c} d\nu  _{\shift (x)}
\leq \frac{M_a}{R^a}\,.
$$
On the other hand, Lemma \ref{3.2 6} applied with $R=R_0$ yields
$$\cL_x^* \nu  _{\shift (x)}  \1 \geq \int _{\overline \D _{R_0}} \cL_x \1 d\nu _{\shift (x)}
\geq  c {R_0}^{-(\al_2 -\tau )t}
\big[\omega^{-1}\big( 8 \log {R_0}\big)\big]^{-\hat \tau  t }\nu _{\shift (x)} (\overline \D _{R_0})\,.
$$
Notice that $\nu _{\shift (x)} (\overline \D _{R_0})\geq \frac12$ since $\nu _{\shift (x)}\in \cM_{\shift (x)}$.
Therefore, 

\begin{align*}
\Phi_x (\nu _{\shift (x)} ) (\overline \D_R^c)=\frac{\cL_x^* \nu  _{\shift (x)} (\overline \D_R^c) }{\cL_x^* \nu  _{\shift (x)}  \1}
\leq &
\frac{2M_a}{c} 
\frac{ {R_0}^{(\al_2 -\tau )t}
\big[\omega^{-1}\big( 8 \log {R_0}\big)\big]^{\hat \tau  t }}{R_0^{a/2}} \frac{R_0^{a/2}}{R^a}
\,.\\
\leq & \frac{2M_a}{c} 
\frac{ \big[\omega^{-1}\big( 8 \log {R_0}\big)\big]^{\hat \tau  t }}{R_0^{b}}
\frac1{R^{a/2}}\,.
\end{align*}
In order to conclude that $\Phi_x(\nu _{\shift (x)})\in \cM_x$ it suffices to set 
$\ep = \frac{a}{2}$ and to show that there exists $R_0$ such that  
$ \frac{2M_a}{c} 
\frac{ \big[\omega^{-1}\big( 8 \log {R_0}\big)\big]^{\hat \tau  t }}{R_0^{b}} \leq 1$.
But this results from an elementary calculation based on the properties of $\omega$:
 $\lim_{r\to \infty} \log (r)/ \log \omega (r) =0$ and $\omega $ is increasing.
\end{proof}

\


\bprop\label{continuity phi}
Let $\cM$ be the invariant set of random measures from Proposition \ref{3.2 4}.
The map $\Phi : \cM \to \cM$ is continuous with respect to the narrow topology.
\eprop

\bpf
Suppose that $\Lambda$ is a directed set and $(\nu^\a)_{\a\in\La}$ is a net in $\cP_m(\J)$ converging to a measure  $\nu \in \cP_m(\J)$ in the narrow topology. If 
$h_{\shift (x),\a} = 1 / \nu^\a _{\shift (x)} (\cL_x\1 )$ then, by \eqref{3.2 1},
$$ 
\Phi_x (\nu^\a_{\shift (x)})= \cL^*_x \left( \frac{1}{\nu^\a_{\shift (x)} (\cL_x\1 )} \nu^n  _{\shift (x)} \right)=
 \cL^*_x \left( h_{\shift (x),\a}\;  \nu^\a_{\shift (x)} \right)
\,.$$
Proposition \ref{25.2 6} implies that $\cL^* $ is continuous with respect to the narrow topology of $\cM(\J )$.
Thus, we have to investigate $h_{\shift (x),\a} \; \nu^\a_{\shift (x)}$. Item (a) of the definition of $\cM$ in
\eqref{3.2 2} and Lemma \ref{3.2 6} imply that there are constants $0<c_1<c_2<\infty$ such that
$$c_1\leq h_{\shift (x),\a} \leq c_2 \quad \text{for all} \quad n\geq 0 \text{ and } x\in X\,.$$
This implies that $(h_{x,\a} \; \nu^\a_{x})_{x\in X}$ is a tight, hence relatively compact, family 
of $\cM_m(\J )$. Let $\mu \in \cM_m(\J )$ be an accumulation point. It is shown (as a matter of fact for sequences but the same argument works for all nets) in 
Lemma~2.9 of \cite{RoyUrb2011} that then 
$$\mu = h \nu$$
for some measurable function $h:X\to (0, \infty )$.
Notice that $\cL^* \mu $ is a random (probability) measure. Hence, the disintegrations of this measure 
$$
\cL^*_x \mu_{\shift (x)} = h_{\shift (x)} \cL^*_x ( \nu _{\shift (x)} )\in \cP (\J_x )
$$
are probability measures. Therefore, $1= \cL^*_x \mu_{\shift (x)} (\1) = h_{\shift (x)} \cL^*_x ( \nu _{\shift (x)} )(\1)$
which implies that the accumulation point $\mu$ is uniquely defined by
$$\mu _{\shift (x)}  = \frac{1}{\cL^*_x ( \nu _{\shift (x)} )(\1)} \nu_{\shift (x)}\quad , \quad x\in X \,.$$
This shows that the net $(\Phi(\nu^\a)_{\a\in\La}$ converges to $\Phi (\nu )$ in the narrow topology. The proof of continuity of $\Phi$ is complete.
\epf

\noindent We are now ready to prove the main result of this section.

\bpf[Proof of Theorem \ref{thm Gibbs}]
Proposition \ref{3.2 4} yields a $\Phi$--invariant convex and compact set $\cM \subset \cP_m(\J )$ of random measures. By Proposition~\ref{continuity phi})
the map $\Phi$ is continuous on $\cM$ for the narrow topology.
Therefore, one can apply Schauder-Tychonoff Fixed Point Theorem in order to get a $\Phi$--invariant random measure
$\nu$. This measure is the required Gibbs state. Finally, the bounds on
$$
\l_x = \nu _{\shift (x)} ( \cL_x\1)
$$
follow again from item (a) of the definition of $\cM$ in
\eqref{3.2 2} and Proposition~\ref{prop bounded operators} together with Lemma~\ref{3.2 6}.
\epf


\

\noindent We have to study these random conformal measures more in greater detail. Here and in the rest of the paper it is very useful to introduce normalized operators
\beq \label{eq nomralized}
\npfx := \l_x^{-1} \pfx,
\eeq
and to employ the notation
$$
\l_x^n = \prod_{j=0}^{n-1} \l_{\shift^j (x)}^n\quad \text{and} \quad
\npfx^n= \l_x^{-n}\pfx^n\,.
$$
We continue to use the radius $R_0$ of the definition of the invariant measure space $\cM$, given in Proposition~\ref{3.2 4}. Clearly we may suppose that $R_0\geq T >0$, $T$ being the constant of Condition \ref{C2}. Condition \ref{C4} is applied to get the following lower estimate.

\blem\label{lem 2.2}
For every $R\geq R_0+1$ and every $0<\d \leq\min \{ \d_0 , 1\}$, where $\delta_0\in (0,1/4)$ comes from Definition~\ref{dfn 1}, there exists $A=A(\d , R, t)\geq 1$ such that 
$$
\nu_x\big(\D(z,\d)\big) \geq A^{-1} \quad  \text{for all} \ x\in X \ \text{and all} \  z\in \julx \  \text{with} \ |z|\leq R.
$$
\elem

\bpf
Covering $\overline\D_{R_0}\cap \julx$ with $\d$--disks whose no more than three elements intersect, and using the fact that $\nu_x(\overline\D_{R_0 })\geq \frac12$, we see that there exists $\tilde a =  a (\d ,R_0)>0$ such that 
\beq \label{eq 3}
\nu_x \big( \D(w_x , \d )\big) 
\geq  \tilde a \quad  \text{for } m-\text{a.e.}\;\;  x\in X\,  \text{ and some }  w_x\in \overline\D_{R_0}
\;\;.
\eeq
Let now $N=N(\d ,R)$ be the number coming from the mixing property of Lemma~\ref{lem: mixing}, applied with $r=\d / 2$ and $R\geq R_0+1$. 
Then, for every $z\in \julx$, $|z|\leq R$,
$$f_x^N \left( \D (z, \d /2 )\right)\supset \overline \D _R\cap \jul_{\shift ^N (x)}\supset \D(w_{\shift^N(x)},\d )\cap \jul_{\shift ^N (x)}\,.$$
Consequently, there exists a holomorphic inverse branch $f_{x,*}^{-N}$ of $f_x^N$ defined on the disk $\D(w_{\shift^N(x)},\d )$
such that $f_{x,*}^{-N}\left(\D(w_{\shift^N(x)},\d )\right) \cap  \D (z, \d /2 )\neq \emptyset $.
We may assume without loss of generality that $N$ is so large that Koebe's distortion theorem together with the expanding property
imply that $\diam \left(f_{x,*}^{-N} (\D(w_{\shift^N(x)}, \d))\right) \leq \d /2 $. Then
$f_{x,*}^{-N}\left(\D(w_{\shift^N(x)},\d )\right) \subset  \D (z, \d  )$. Hence, 
\begin{equation*}
\begin{split}
\nu_x(\D(z, \d ) ) 
& \geq \nu_x\left(f_{x,*}^{-N}\left(\D(w_{\shift^N(x)},\d )\right)\right)\\
&\geq \frac1K \l_x^{-N} \left| (f_x^N)'(f_x^{-N} (w_{\shift^N(x)}) \right|_\tau ^{-t} \nu_{\shift^N(x)}(\D(w_{\shift^N(x)},\d )) \\
&\geq a 
\end{split}
\end{equation*}
for some $a=a(t,R, \d ) >0$ by \eqref{eq 3}, Condition \ref{C4}, and since $\l_x\leq C<\infty$ for all $x\in X$.
\epf

\smallskip

\section{Uniform bounds and invariant densities}

We are now able to prove the following uniform bound for the normalized operators $\npfx^n$.

\bprop \label{prop1}
There exists $M=M_t<\infty$ such that 
$$
\|\npfx^n\|\leq M\quad \text{for every} \ n\geq 1 \text{ and } m-\text{a.e. } x\in X \,.
$$
\eprop

\

\noindent
By combining this result with Proposition \ref{prop two norm}, we obtain the following:

\bcor\label{uniform two norm} For every $x\in X$, every $g\in \cH_\b (\julx )$, and all $n\geq 1$, we have
\beq \label{l21.20}
v_\b (\npf_x ^ng)\leq M\left( \|g\|_\infty + K (c\g^n)^{-\b}v_\b(g)\right).
\eeq
\ecor

\noindent Proposition \ref{prop1} and Corollary \ref{uniform two norm} together imply  that the above uniform bound 
is also valid with respect to the H\"older norm $\|\cdot\|_\b$. For simplicity we will use the same bound $M$ in the sequel.

\bcor \label{uniform holder norm}
There exists $M=M_t<\infty$ such that 
$$\|\npfx^n\|_\b\leq M\quad \text{for every} \quad n\geq 1 \ \text{ and }\  m-\text{a.e. } x\in X \,.$$
\ecor
\

\noindent
We first need an auxiliary result. Let $0<\d \leq \d_0$ and $R\geq R_0+1$, where, we recall, $\delta_0\in (0,1/4)$ comes from Definition~\ref{dfn 1} and $R_0>0$ is taken from Proposition~\ref{3.2 4}.

\blem \label{l21.1}
For every $n\geq 1$ we have 
$$
\npfx^n \1 (w)\leq K A \quad , \quad w\in \jul_{\shift ^n(x)} \;\; , \;\; |w|\leq R\; , \;\; x\in X\, ,$$
where K is the distortion constant from Lemma \ref{lem dist operator} and A is defined in Lemma \ref{lem 2.2}.
\elem

\bpf
Let $w\in \jul_{\shift ^n(x)}$, $ |w|\leq R$. Then, using Lemma \ref{lem dist operator} and Lemma \ref{lem 2.2}, we get
$$
1= \int \npfx^n \1 d\nu_{\shift ^n(x)} \geq  \int_{\D(w,\d )} \npfx^n \1 d\nu_{\shift ^n(x)}
\geq \frac1K \npfx^n\1 (w) \nu_{\shift ^n(x)}\big( \D(w,\d )\big)\geq \frac{\npfx^n\1 (w)}{AK}\,.
$$
\epf

\bpf[Proof of Proposition \ref{prop1}]
From the lower bound on $\l_x$ given in Theorem \ref{thm Gibbs} and from the fact that $\pfx \1(w)\to 0$
as $|w|\to\infty$ uniformly in $x\in X$ (see Proposition \ref{prop bounded operators}), it follows that there exists $R\geq R_0+1$ such that 
\beq\label{l21.11}
\npfx\1 (w) \leq 1 \quad \text{for} \quad w\in \jul_{\shift (x)}\cap \D_R^c\,
\eeq
for $m$-a.e. $x\in X$.
\begin{claim}
Set $M=K A(\d_0 , R,t )$, again with constants as in Lemma~\ref{l21.1}. Then $\npf_x^n\1 \leq M$ for every $n\geq 1$ and $m$-a.e. $x\in X$.
\end{claim}
\noindent It suffices to prove this claim. It will be done by induction. In the case $n=1$ the results follow directly from \eqref{l21.11} and Lemma \ref{l21.1}. So, fix
$n\geq 1$ and suppose that Claim holds for this $n$. We have to show that
$$
\npf^{n+1}_{\shift ^{-(n+1)}(x)} \1(w)\leq M\quad \text{for every}\quad w\in \julx \;\; \text{and } a.e. \;x\in X\,.
$$
If $w\in \jul_{x}\cap \D_R$, then it suffices to apply Lemma  \ref{l21.1}.
Otherwise, i.e. if $|w|\geq R$, then 
$$\npf^{n+1}_{\shift ^{-(n+1)}(x)}\1 (w) = \npf_{\shift^{-1} (x)}\left(\npf^{n}_{\shift ^{-(n+1)}(x)}\1\right) (w)
\leq M  \npf_{\shift^{-1} (x)} \1 (w) \leq M\,.$$
\epf

%
%
\section{Invariant positive cones and Bowen's contraction}

G. Birkhoff in \cite{Bir57} reinterpreted Hilbert's pseudo-distance on positive cones in a way which allowed him to
show that every linear map preserving cones is a weak contraction. This enabled him to
give an elegant proof of the Perron-Frobenius theorem based on Banach's contraction principle.
Various versions of Ruelle's Perron-Frobenius theorem have been obtained since then using Birkhoff's strategy
(see, for example, 
Liverani \cite{Liv95} and Rugh \cite{Rug08, Rug10} who, at the same time, considered random dynamics and introduced a complexification scheme leading to real analyticity of the dimension).

In our setting, with unbounded phase spaces $\julx$, we encounter several problems.
First of all, because of the behavior of the functions at infinity, every reasonable invariant cone contains many functions that all are at the infinite Hilbert distance from each other. These cones have many, in fact uncountably many, connected components
that are at finite distances from each other. The second problem is that it is hard to get a strict-contraction property since
the mixing property which is at our disposal (Lemma \ref{lem: mixing}) is too weak; one only has mixing on bounded regions. 

Our way to overcome these difficulties is to define appropriate invariant cones and then to avoid Birkhoff's strategy, but instead, to employ an argument inspired by Bowen's lemma \cite[Lemma 1.9]{Bow75}.
For compact phase spaces this lemma is indeed equivalent to a strict contraction in the Hilbert metric.
In our situation this is not the case but it turns out that Bowen's lemma is sufficiently tricky so that we can use
some appropriate version of it that leads to the following exponential convergence result. 

From now on the number $\d_0>0$ in the definition of the variation of H\"older functions will be replaced 
by a smaller number $0<\d\leq \d_0$ as explained in \eqref{l21.8}.

\bthm\label{thm expo} Let $(f_x)_x$ be a hyperbolic transcendental random system. We then have the following.
\ben
\item There exists a unique $\rho \in \cH_\b (\J ) $, which is an invariant density, i.e. $\npf \rho = \rho$.
\item There are $B>0$ and $\vth \in (0,1)$ such that
$$ 
\left\|\npfx ^n g_x -\int g_x \,d\nu_x \, \den_{\shift^n(x)}\right\|_\b
\leq B\vth ^n \| g_x\|_\b
$$
for every $ g_x \in H_\b (\julx )$ and a.e. $ x\in X.$ 
\een
\ethm

\brem 
Multiplying, as usually, the random Gibbs state $\nu$ by the invariant function $\rho $ of Theorem \ref{thm expo},  (1)
gives again a unique, by  Theorem \ref{thm expo} (2), invariant random Gibbs state $\mu\in \cP_m(\J )$ whose disintegrations are
\beq \label{equilibrium states}
\mu_x = \rho_x \nu_x \quad , \quad x\in X.\eeq
Moreover,  $\mu$ is ergodic. Indeed, if there existed an invariant set $E\subset \J$ with $0<\mu (E) <1$, then
$\mu_1 = \1_E \mu$ and $\mu_2= \1_{E^c} \mu$ would be two invariant random Gibbs states. But this would contradict the above uniqueness property.
\erem

\brem
Notice that, as a straightforward consequence of the assertion (2) of this theorem, we also get
 exponential convergence for random H\"older observables, i.e in
$  \cH_\b (\J ) $, with respect to the canonical norm of this space:
$$ | \npf ^n g - \pi _\rho ( \npf ^n g )|_\b \leq  B\vth ^n  |g|_\b \quad , \quad g\in  \cH_\b (\J )\, , $$
where $\pi_\rho : \cH_\b (\J ) \to \; <\rho>$ is the canonical projection defined by $\pi_\rho (g)_x = \int g_x \, d\nu_x \rho_x$.

In particular, $\npf ^n\1 \longrightarrow \rho $ exponentially fast.
\erem


\subsection{Invariant cones}
Consider the following cones:
\beq\label{l21.5}
\cC_x := \left\{g\geq 0 \; :\;\;  \|g\|_\infty \leq \cA \int g d\nu_x  <\infty\; \; \text{and} \;\;
 v_\b (g)\leq H \int g d\nu_x\right\}\,.
\eeq
\beq\label{l21.6}
\cC_{x,0} := \left\{g\in \cC_x \; :\;\; g\leq 2 M_t\cA \left(\int g d\nu_x\right) \; \npf _{\shift^{-1} (x)}\1 \right\}\,.
\eeq
Since we are primarily interested in the projective features of these cones, it is convenient for us to use  the following sections
\beq\label{l21.6slice}
\La _x = \{ g\in \cC_x \; , \;\; \nu_x (g) =1\} \quad \text{and} \quad \La _{x,0} =\La_x\cap \cC_{x,0} 
\; \; , \;\; x\in X\,.
\eeq
Hence both type of cones do depend on constants $\b \geq 1$, $\cA>0$, $H>0$, and even, indirectly, on $\d$ which was defined in the paragraph between Corrolary~\ref{uniform holder norm} and Lemma~\ref{l21.1}. This dependece comes via the value of H\"older variation $v_\b$. In fact, when we deal in the sequel with H\"older functions $g$, then we assume that the variation $v_\b (g)$ is evaluated on disks of radius $\d$, i.e.  $\d_0$ is replaced by $\d$ in \eqref{l21.9}. Whenever the dependence on the constants is important we will indicate this and write $\cC_x (\cA, H)$  or even  $\cC_x (\cA, H, \b , \d )$, and similarly for the second type of cones.
In order to produce cones with good properties, for example invariance, we have to choose carefully these constants.

We continue to write $M=M_t$ for the uniform bound given in Proposition \ref{prop1} and in Corollary \ref{uniform holder norm} and $K=K_t$ for the distortion constant appearing in Lemma \ref{lem koebe} and Lemma \ref{lem dist operator}.
First of all, let $0<\d \leq \d_0$ be such that
\beq \label{l21.8}
\frac12 +\left(2MK+4 \right)\d^\b \leq 1\,.
\eeq
The radius $R_0$ has been defined in Lemma \ref{3.2 6}. Increasing it if necessary we may suppose
that Lemma \ref{lem 2.2} is valid with $R=R_0$ and, for the same reason as in \eqref{l21.11}, that
\beq\label{l21.10}
2M \npfx \1 \leq 1 \quad \text{in} \quad  \D_{R_0}^c \cap \julx \; , \;\; x\in X\,.
\eeq
Define now, with $A(\d , R , t)$ from Lemma \ref{lem 2.2},
\beq \label{l21.13}
\cA:= 2\max\{ 1, A(\d , R_0 , t),M \} \quad \text{and} \quad H= 2MK\cA+4\,.
\eeq
Notice that $\cA \geq 1$. This ensures that the constant function $\1 \in \cC_x$, $x\in X$. Finally, let $N_0\geq 1$ be such that 
\beq \label{m22.08}
MK(c\g^{N_0})^{-\b } H \leq 1\,.
\eeq

\bprop \label{l21.7}
With the above choice of constants and  for every $n\geq N_0$,
$$
\npfx^n \left( \cC_x\right) \subset \cC_{\shift ^n (x),0} \subset \cC_{\shift ^n (x)} \; , \;\; x\in X\,.
$$
\eprop

\bpf
Let $g\in \cC_x$. We may assume that $\int g d\nu_x =1$.  
We will show that $\npfx^n g \in \cC_{\shift ^n (x) , 0}$ for every $n\geq N_0$.
Let in the following $n\geq N_0$. From the two-norm type inequality \eqref{l21.20} and from the definition of the cone, we get that
\beq \label{l21.21}
v_\b (\npf_x^n g)\leq M\left( \cA + K (c\g^n)^{-\b}H\right)\leq M\cA+1\leq H,
\eeq
where the last two inequalities result from the choice of $N_0$ and from the definition of $H$.
Then,
\beq \label{l21.22}
\npfx^n g = \npf_{\shift^{n-1}(x)}\left( \npfx ^{n-1} g\right)\leq M\|g\|_\infty  \npf_{\shift^{n-1}(x)}\1
\leq M\cA \npf_{\shift^{n-1}(x)}\1\,.
\eeq
In order to see that $\pfx^ng\in \cC_{\shift ^n (x),0} $ it remains to estimate $\| \npfx^n g\|_\infty$.
Since we already have \eqref{l21.21}, we obtain, for every $|z|\leq R_0$, the following:
\begin{equation*}\begin{split}
1&=\int \npfx ^ng \, d\nu_{\shift^n (x)}
\geq \int_{\D(z,\d )}\npfx^n g\, d\nu_{\shift^n (x)} \\
&\geq \left( \npfx^n g(z) -H\d^\b \right)\nu_{\shift^n (x)} (\D (z,\d ))
\geq \left( \npfx^n g(z) -H\d^\b \right) A(\d , R_0,t)^{-1} \,,
\end{split}\end{equation*}
where the last inequality holds true due to Lemma~\ref{lem 2.2}. Therefore, 
$$ 
\npfx^n g(z)
\leq \cA /2+H\d^\b\leq  \cA \left(\frac12 + (2MK+4)\d^\b \right)
\leq \cA,
$$
by \eqref{l21.8}. If $|z|\geq R_0$, then  it suffices to combine \eqref{l21.22} and \eqref{l21.10}
in order to conclude this proof with the inequality,
$$
\npfx^ng(z)\leq M\cA \npf_{\shift^{n-1}(x)}\1 (z)\leq \cA\,.
$$
The proof is complete.
\epf

\subsection{Cone contraction via Bowen's lemma.}

Let $R_1\geq R_0$ be such that
\beq\label{m22.01}
2AM \npfx\1 \leq 1  \quad \text{in} \quad \D_{R_1}^c\,.
\eeq

\blem\label{m22.02}
For every $R\geq R_1$ there are $N=N_R\geq N_0$ and $a=a_R>0$ such that 
$$
\npfx^N g_{\big| \D_{2R}}\geq a \quad \text{for every} \quad g\in \La_{x,0}\;, \;\; x\in X \,.
$$
\elem

\bpf
Let $g\in \Lambda_{x,0}$. Since $\int gd\nu_x =1$, we have that $\|g\|_\infty \geq 1$. Hence, by the choice of $R_1$,
$$
g\leq 2M\cA \, \npfx\1 \leq 1 \quad \text{in} \quad \D_{R_1}^c.
$$
Thus, there exists $z_{max}\in \overline \D_{R_1}$ with $g(z_{max})=\|g\|_\infty \geq 1$.

Let $0<r\leq \d$ be such that $Hr^\b \leq \frac14$. The mixing property of Lemma \ref{lem: mixing}
implies the existence of $N= N(r,R)\ge 0$ such that every $w\in \jul_{\shift^N(x)}\cap  \D_{2R}$
has a preimage $z_0\in f_x^{-N} (w) \cap \D(z_{max}, r)$. Therefore, using  Condition \ref{C4}, for every such $w$, we get that
\begin{equation*}
\begin{split}
\npf_x^N g(w) 
&\geq \l _x^{-N}\left|(f_x^N)'(z_0) \right|_\tau^{-t}g(z_0) \geq C^{-N} \left|(f_x^N)'(z_0) \right|_\tau^{-t}
\left( g(z_{max})-Hr^\b\right)\\
&\geq C^{-N} \inf _{
|z|\leq R_1\; |f_x^N (z)|\leq 2R
} \left|(f_x^N)'(z) \right|_\tau^{-t} \left(1-\frac14\right) \\
&=: a>0.
\end{split}
\end{equation*}
The proof is complete.
\epf

Notice that there is no way to get a global, valid on the whole Julia set, version of Lemma \ref{m22.02}.
This is why we have to work with the following truncated functions. We remark that our cones are chosen in such a way that such truncations can be made to lie inside them. This is not the case for the standard Bowen cones.

Let $\ph_1: \C\to [0,1]$ be a Lipschitz function such that $\ph_1\equiv 1$ on $\D_1$ and  
$\ph_1\equiv 0$ on $\D_2^c$. For $R\geq 1$ define $\ph_R(z)=\ph_1(z/R)$. Then $\ph_R$ is also Lipschitz with variation $v_1(\ph_R)\to 0$ as $R\to \infty$. Define 
\beq\label{m22.04}
\ph_{x,R} :=  \ph_R \npf_{\shift^{-1}(x)}\1  \,.
\eeq
Then, $0\leq \ph_{x,R} \leq \npf_{\shift^{-1}(x)}\1$, $\ph_{x,R} \equiv \npf_{\shift^{-1}(x)}\1$ on $\D_R$ and $\ph_{x,R} \equiv 0$
in $D_{2R}^c$. The functions $\ph_{x,R}$ are Lipschitz with $v_1(\ph_{x,R} )\to v_1 (\npf_{\shift^{-1}(x)}\1)$
uniformly as $R\to \infty$. Therefore, given the definition of the cones, especially the definitions of the constants $\cA , H$ in \eqref{l21.13}, and the formulas established in the course of the proof of Proposition~\ref{l21.7},
it follows that 
$$\ph_{x,R} \in \cC_{x,0}\quad \text{, $x\in X$,}$$
provided that $R$ is sufficiently large. We will assume that $R_1$
is chosen so that these truncated functions belong to the cones for all $R\geq R_1$.
Suppose also, in what follows, that $\eta>0$ is chosen such that
\beq\label{m22.03}
0<\eta\leq \min \left\{ \frac13 , \; \frac1H , \; \frac12\frac{a}{M}\right\} \,.
\eeq
With these choices we will now obtain the following version of
Bowen's result \cite[Lemma 1.9]{Bow75}.

\blem\label{m22.04}
For every $R\geq R_1$ and with $N=N_R\geq N_0$ given by Lemma \ref{m22.02},
$$ \npfx^N g - \eta \ph_{{\shift^N(x)},R} \in \cC_{{\shift^N(x)},0} \quad \text{for every}\quad g\in \La_{x,0}\,.$$
\elem

\bpf
Let $x\in X$, let $g\in \La_{x,0}$, and let $R\geq R_1$.
Lemma \ref{m22.02} shows that for $0<\eta <\frac12\frac{a}{M}$,
$$
\npfx^N g -\eta \ph_{\shift^N(x)} >\frac{a}{2} >0 \quad
\text{on} \quad \D_{2R}\cap \jul_{\shift^N(x)} \,.
$$
Set
\beq\label{m22.05}
g'=\frac{\npf_x^N g -\eta \ph_{{\shift^N(x),R}}}{1-\eta_{{\shift^N(x)},R}} \;\; \;\text{where}\;\; \;
\eta_{{\shift^N(x)},R}:= \eta \int\ph_{{\shift^N(x)},R}d\nu_{\shift^N(x)} \,.\eeq
Then $\int g' d\nu_{\shift^N(x)}  =1$ and $g'>0$. We have,
by
$$(1-\eta_{{\shift^N(x)},R})\,g'\leq  M\|g\|_\infty \npf_{\shift^{(N-1)}(x)}\1 + \eta  \npf_{\shift^{(N-1)}(x)}\1
\leq (M\cA +\eta )   \npf_{\shift^{(N-1)}(x)}\1 \,.$$
But $0<\eta_{{\shift^N(x)},R}\leq \eta\leq \frac13$ and thus
$g'\leq 2M\cA   \npf_{\shift^{(N-1)}(x)}\1 $. This means that the function $g'\in \La_{{\shift^N(x)},0} $ provided that we
can show that $g'\in \La_{{\shift^N(x)}} $.

In order to estimate the variation of $g'$ we use again the two-norm type inequality \eqref{l21.20}:
$$v_b(g') \leq \frac{1}{1-\eta_{{\shift^N(x)},R}} \left(M\|g\|_\infty  +MK (c\g ^N)^{-\b}v_\b (g)+ \eta v_\b  (\ph_{_{\shift^N(x)},R})\right)\,.
$$
Remember that $g,\ph_{_{\shift^N(x)},R}\in \cC_x$, that $\eta \leq \min \left\{ \frac13 , \; \frac1H \right\}$, and that we have \eqref{m22.08}. Therefore,
$$v_b(g') \leq 2 \left(M\cA  +1+ 1\right)=2M\cA+4\leq H\,.
$$
It remains to estimate $\|g'\|_\infty$. If $z\in \jul_{\shift^N(x)} \cap \D_{R_0}$, then
$$
1=\int g' d\nu_{\shift^N(x)}
\geq \int _{\D(z, \d)}g' d\nu_{\shift^N(x)} 
\geq (g'(z)- H\d^\b ) \nu _{\shift^N(x)} (\D(z,\d ))\,.
$$
Using once more Lemma \ref{lem 2.2} and the choice of $\d $ in \eqref{l21.8}, we obtain
$$ 
g' (z)
\leq A(\d , R_0 , t) + H\d^\b 
\leq \frac{\cA}{2} + (2MK\cA+4)\d^\b 
\leq \cA\left(\frac12 +(2MK+4)\d^\b \right) 
\leq \cA\,.
$$
If $z\in \jul_{\shift^N(x)} \cap \D_{R_0}^c$, then $g'(z)\leq 2M\cA   \npf_{\shift^{(N-1)}(x)}\1 (z)\leq \cA$
by the choice of $R_0$ (see \eqref{l21.10}).
The proof is complete.
\epf

Applying repeatedly Lemma  \ref{m22.04}  gives the desired contraction.

\bprop\label{m22.09}
For every $\ep >0$ there exists $n_\ep \geq 1$ such that for every $n\geq n_\ep$ and
a.e. $x\in X$,
\beq\label{28.2 1}
\left\|\npf _x^{n} g_x -\npf _x^{n} h_x\right\|_\b \leq\ep  \quad \text{for all} 
\quad g_x,h_x\in \La_{x,0} \,.\eeq
\eprop

\bpf
Let $R\geq R_1$ and $N=N_R\geq N_0$ be like in Lemma \ref{m22.04}.
and let $g=g_x^{(0)}\in \La _{x,0}$. With the notation of the previous proof, and in particular with
the numbers $\eta_{{\shift^N(x)},R}$ defined in \eqref{m22.05}, we get from
Lemma \ref{m22.04} that
$$\npfx ^N g = \eta \ph_{{\shift^N(x)},R} + (1- \eta_{{\shift^N(x)},R})g_{\shift ^N (x)}^{(1)}$$
for some $g_{\shift ^N (x)}^{(1)}\in \La _{\shift ^N (x),0}$. Applying $\npf_{\shift ^N(x)} ^N$ to this equation
and using once more Lemma \ref{m22.04} gives
$$\npfx ^{2N} g=  \eta \npf_{\shift ^N(x)} ^N\ph_{{\shift^N(x)},R}
+ (1- \eta_{{\shift^N(x)},R})\eta \ph_{{\shift^{2N}(x)},R}
+ (1- \eta_{{\shift^N(x)},R})(1- \eta_{{\shift^{2N}(x)},R}) g_{{\shift^{2N}(x)}}^{(2)}$$
for some $g_{\shift ^{2N} (x)}^{(2)}\in \La _{\shift ^{2N} (x),0}$. Inductively it follows that for every $k\geq 1$
there is a function $g_{\shift ^{kN} (x)}^{(k)}\in \La _{\shift ^{kN} (x),0}$ such that
$$\npfx ^{kN} g= \eta \sum_{j=1}^k \left( \prod_{i=1}^{j-1} (1-  \eta_{{\shift^{iN}(x)},R})\right)
\npf_{\shift^{jN}(x)}^{(k-j)N}\ph_{{\shift^{jN}(x)},R} +  \prod_{i=1}^{k} (1-  \eta_{{\shift^{iN}(x)},R})\,
g_{\shift ^{kN} (x)}^{(k)}\,.
$$
Observe that the first of these two terms does not depend on $g$. Therefore, for every $g,h\in \La_{x,0}$
there are $g_{\shift ^{kN} (x)}^{(k)}, h_{\shift ^{kN} (x)}^{(k)}\in \La _{\shift ^{kN} (x),0}$ such that
\begin{equation}\label{120140501}
\npfx ^{kN} g -\npfx ^{kN} h =  \prod_{i=1}^{k} (1-  \eta_{{\shift^{iN}(x)},R})\,
\left( g_{\shift ^{kN} (x)}^{(k)} - h_{\shift ^{kN} (x)}^{(k)}\right)\,.
\end{equation}
Remember that $\eta _{y,R} = \eta \int \ph _{y,R} d\nu _y \geq \eta \int _{\overline \D_{R_0}} \npf_{\shift^{-1} (y)}\1 d\nu _y$ for all $R\geq R_0$, and that, by
Lemma \ref{3.2 6}, there exists a constant $c=c(R_0)>0$ such that $\npf_{\shift^{-1} (y)}\1 \geq c$ on $\overline \D_{R_0}$. Therefore,
$$
1 > \eta 
\geq \eta _{y,R} 
\geq \eta c \nu _y (\overline \D_{R_0}) 
\geq \eta \frac{c}{2} 
= \tilde \eta>0.
$$
Thus,
$$
1-\eta _{y,R}\le 1-\tilde \eta.
$$
Along with \eqref{120140501}, this allows us to deduce the the uniform bound of Proposition~\ref{prop1}, with some $n_{1, \ep}\ge 1$ sufficiently large, for the supremum norm rather than the H\"older one.
In order to get the appropriate estimate for the $\b$--variation we need once more  \eqref{l21.20}.
Write $n=m+n_{2,\ep}+n_{1, \ep}$ with some $n_{2, \ep}$ to be determined in a moment and some $m\ge 0$. Then for all $g,h\in \La_{x,0}$, we have
$$
\begin{aligned}
v_\b \left(\npf_x^n g -\npf_x^nh \right)
&= v_\b \left( \npf_{\shift ^{n_{1, \ep}}(x)}^{m+n_{2, \ep}} \left( \npfx ^{n_{1, \ep}}( g -h)\right)\right)\\
& \leq M\left( \left \| \npfx ^{n_{1, \ep}}( g -h)\right\| _\infty  + K(c\g ^{m+n_{2, \ep}} )^{-\b}v_\b \left( \npfx ^{n_{1,\ep}}( g -h)\right)\right) \\
&\leq  M\ep + MK (c\g ^{n_{2, \ep}})^{-\b} \, 2H,
\end{aligned}
$$
since $\npfx ^{n_{1, \ep}} g, \npfx ^{n_{1, \ep}} h \in \cC_{\shift ^{n_{1, \ep}}(x) , 0}$. It suffices now to choose $n_{2, \ep}\ge 0$ sufficiently large in order to conclude this proof.
\epf

\

\bpf[Proof of Theorem \ref{thm expo}  (1)] 
Consider $\rho ^k = \npf ^{k} \1$. First of all, Proposition \ref{l21.7} implies that $\rho ^k_x\in\La _{x,0}$
for every $k\geq N_0$. Hence Proposition \ref{m22.09} applies and gives
$$\|\rho ^k_x - \rho ^l_x \|_\b \leq \ep \quad \text{for every} \quad 
l\geq k \geq n_\ep \;\; , \;\; x\in X\,.$$
This shows that $(\rho ^k_x)_k$ is a, uniformly in $x\in X$, Cauchy sequence of $(\cH_\b (\J_x ) , \| .\|_\b)$
and hence there is a limit $\rho \in \cH_\b (\J)$. Clearly, $\npf \rho = \rho$ and $\rho_x\in \La _{x,0}$, $x\in X$.
Uniqueness of this function follows from the contraction given in \eqref{28.2 1}.
\epf

\sp\bpf[Proof of Theorem \ref{thm expo} (2)]
Since $ \cA , H \ge 2$, we have that
$$
\big\{\1+h_x: \;\|h_x\|_\b < 1/4 \big\} \subset \cC_x,
$$
for all $x\in X$.
Let $g\in H_\b(\julx )$, $g\not\equiv 0$ be arbitrary. Then
$$
h:= \frac{g}{8\|g\|_\b} =(h+\1)- \1
$$
is a difference of functions from $\cC_x$. If $\ep >0$ and $n=n_\ep $ is given by Proposition \ref{m22.09}, then
\begin{equation*}\begin{split}
\Bigg\| \npfx^n \Bigg(h -&\left(\int h d\nu_x\right) \den_{\shift^n(x)}\Bigg) \Bigg\|_\b
\le \left\| \npfx^n h -\left(\int h d\nu_x\right) \den_{\shift^n(x)} \right\|_\b \le \\
&\leq \left\| \npfx ^n (\1+h) - \int (\1+h) \, d\nu_x \, \den_{\shift^n(x)} \right\|_\b +\left\| \npfx^n\1 -\den_{\shift^n(x)} \right\|_\b\\
&\leq \ep  \int (\1+h) \, d\nu_x + \ep \\
&\leq \frac{17}8\ep.
 \end{split}\end{equation*}
This shows that for every $\ep >0$ there exists $N=N_\ep$ such that
$$
\left\| \npfx^N \left(g -\left(\int g d\nu_x\right) \den_{\shift^N(x)}\right) \right\|_\b \leq \ep \| g\| _\b \quad \text{for every} \quad g\in H_\b (\julx)\, .$$
Fix $\ep:=1/2$ and let $N=N_{1/2}$. Write any integer $n\ge 0$ in a unique form as $n=kN+m$, where $k\ge 0$ and $m\in \{0,...,N-1\}$.  Then, for every $g\in H_\b (\julx )$ 
we have,
\begin{equation*}\begin{split}
\left\| \npfx^n g -  \int g d\nu_x \den_{\shift^n(x)} \right\|_\b &= 
\left\| \npf_{\shift^{kN}(x)}^m \left( \npfx^{kN}\big( g -\int g d\nu_x \den_x \big) \right)\right\|_\b\\
&\leq M \left(\frac12\right)^k \left\| g -\int g d\nu_x \den_x\right\|_\b\\
&\leq 2 M  \left(\frac1{2^{1/N}}\right)^n \left( 1+\|\den_x\|_\b\right) \| g\|_\b \,.
 \end{split}\end{equation*}
This completes the proof of Theorem \ref{thm expo}.
\epf

%
%

\section{Exponential decay of correlations and CLT}

Exponential decay of correlations is now a fairly straightforward consequence of Theorem~\ref{thm expo} (2).
It will be valid for functions of the following spaces. 

Let 
$\cH_\b^p (\J)$ be the space of functions $g:\J\to \R$ with H\"older fibers $g_x\in \cH_\b (\J_x )$ and such that
$\|g_x\|_\b\in L^p (m)$. The canonical norm is 
$$|g|_{\b , p}= \left( \int_X  \|g_x\|_\b ^p dm(x)\right)^\frac{1}{p}\,.$$
Replacing in this definition the $\b$--H\"older condition on the fiber $\J_x$ by a $L^1 (\nu_x)$ condition leads to a space of functions that will be denoted by $L_{\nu }^{ 1, p}(\J)$. The natural norm is in this case 
$$| g|_{\nu }^{ 1, p} = \left(\int _X \| g_x\|^p _{L^{1} (\nu_x )} dm(x)\right)^\frac1p \,.$$
Clearly, if $p=1$ then $L_{\nu }^{ 1, 1}(\J)=L^1 (\nu)$. 
In both cases we also consider $p=\infty$ and then the $L^p$ norms are replaced by the sup--norm. 

\bthm\label{thm-correlations}
Let $(f_x)_x$ be a hyperbolic transcendental random system and let $p,q\in [1,\infty ]$ such that $\frac1p+\frac1q=1$.
Then, for every  $g\in L_{\nu }^{ 1, p}(\J)$, $h\in \cH_\b ^q (\J)$ with $\int_{\cJ_x}h_x\,d\mu_x=0$ and for every $n\geq 1$, we have
\begin{align*}
\left| \int_{\cJ}(g\circ f^n) \, h \,d\mu \right| 
= \left| \int _X \int _{\J_x}(g_{\shift^n(x) }\circ f_x^n) \, h_x d\mu_x dm(x) \right|
\leq b \vartheta ^n \; |g|_{\nu}^{ 1, p} \;  |h|_{\b,q}
\end{align*}
for some positive constant $b$ and some $\vartheta\in(0,1)$. 
\ethm

\begin{proof}
A standard calculation and application of Theorem~\ref{thm expo} (2) gives
$$
\begin{aligned}
\left|\int_{\J_x}  (g_{\shift^n(x) }\circ f_x^n) \, h_x d\mu_x  \right| &= 
\left|\int_{\J_{\shift ^n (x)}}g_{\shift^n(x) } \, \npf _x^n \left( h_x\den_x \right) d\nu_{\shift^n (x)}  \right| \\
&\leq \left\| \npf _x^n \left( h_x\den_x \right) \right\|_\b \| g_{\shift^n(x) } \|_{L^1 (\nu_{\shift^n (x)} )} \\
&\leq   B \vartheta ^n \left\|  h_x\den_x\right\|_\b \| g_{\shift^n(x) } \|_{L^1 (\nu_{\shift^n (x)} )}
\leq b \vartheta ^n \left\|  h_x\right\|_\b \| g_{\shift^n(x) } \|_{L^1 (\nu_{\shift^n (x)} )}
\end{aligned}
$$
for some constant $b>0$ since $\|\den _x\|_\b\leq M$ for all $x\in X$ by Corollary \ref{uniform holder norm} and Theorem~\ref{thm expo} (2). Therefore,
$$
\begin{aligned}
\left| \int_{\cJ}(g\circ f^n) \, h \,d\mu \right| 
&= \left| \int _X \int_{\J_x}(g_{\shift^n(x) }\circ f_x^n) \, h_x d\mu_x dm(x) \right|\\
&\leq b \vartheta ^n  \int _X   \left\| h_x \right\|_\b \| g_{\shift^n(x) } \|_{L^1 (\nu_{\shift^n (x)} )}  dm(x) \\
&\leq b \vartheta ^n \left(  \int _X   \left\| h_x \right\|_\b ^q  dm(x) \right)^\frac1q \left(  \int _X   \| g_{x} \|_{L^1 (\nu_{x} )}^p  dm(x) \right)^\frac1p\, .
\end{aligned}
$$
\end{proof}

Finally, following Gordin and Liverani's method, one can obtain various versions of the central limit theorem (CLT).
Here is the simplest one.

\bthm\label{CLT}
Let $\psi \in \cH_\b (\J) \cap L^\infty (\J )$ such that
$ \int _{\J_x} \psi _x d\mu_x =0$, $x\in X\,.$
If $\psi$ is not a coboundary in $L^2(\J,\mu)$ (meaning that there is no $u\in L^2(\J,\mu)$ such that $\psi=u-u\circ F$), then there exists $\sg>0$ such that, for every $t\in \R$, 
$$
\mu \left( \{z\in \J \, ; \; \frac1{\sqrt{n}} S_n \psi (z) \leq t\}\right) \to 
\frac1{\sg \sqrt{2\pi }}\int _{-\infty}^t exp (-u^2/2\sg^2 ) \, du
$$
\ethm

\bpf
The dual operator 
$U_x^*:L^2 (\J_x , \mu _x ) \to   L^2 (\J_{\shift (x)},\mu _{\shift (x)})$ of the Koopman operator
$U_x\psi_x = \psi_x \circ f_x$ is  given by
$$U_x^* \psi_x= \frac1{\rho_\shift (x)}\npf _x (\den_x \psi_x )\,.$$
By Gordin's result \cite{Gordin69} it suffices to check that $\sum_k \|U^k U^{*k}\psi\| _{L^2(\mu )}<\infty $.
We have
$$
\| U^k U^{*k}\psi\|^2_{L^2 (\mu )}= \int _\J (U^{*k}\psi)^2 \circ f^k\, d\mu = \int _\J (U^{*k}\psi)^2\, d\mu
$$
by invariance of the measure $\mu$. Therefore,
$$
\| U^k U^{*k}\psi\|^2_{L^2 (\mu )}=  \int _\J \psi \, U^kU^{*k}\psi\, d\mu
\leq \|\psi \| _\infty  \int _\J |U^{*k}\psi | \circ f^k\, d\mu=\|\psi \| _\infty  \int _\J |U^{*k}\psi | \, d\mu
$$
by the same argument. Now,
$$
\begin{aligned}
 \int _\J |U^{*k}\psi | \, d\mu =&
  \int _X \int _{\J_x}  \frac1{\rho_x}\left|\npf_{\shift^{-k}(x)}(\rho_{\shift^{-k}(x)}\psi_{\shift^{-k}(x)})\right|d\mu_x dm\\
  =&  \int _X \int _{\J_x} \left|\npf_{\shift^{-k}(x)}(\rho_{\shift^{-k}(x)}\psi_{\shift^{-k}(x)})\right|d\nu_x dm\,.
\end{aligned}
$$
The fibers of $\psi$ having $\mu_x$--integral zero, it follows from Theorem~\ref{thm expo} (2) as in the preceding proof that, for some constant $b>0$, 
$$ \int _\J |U^{*k}\psi | \, d\mu\leq b \vth^k\int _X \|\psi_{\shift^{-k}(x)}\|_\b dm=  b \vth^k |\psi | _\b\,.
$$
In conclusion,
$$
\| U^k U^{*k}\psi\|^2_{L^2 (\mu )}\leq  b \vth^k |\psi | _\b \| \psi \| _\infty
$$
which directly implies Gordin's $L^2$--summability condition.
\epf

\

%
%

\section{Appendix: Facts from Nevanlinna Theory and uniform bounds of transfer operators} 
\label{sec: nev}

\subsection{FMT and proof of Proposition \ref{prop bounded operators}}

The goal here is to establish the uniform bounds of the transfer operators 
claimed in Proposition \ref{prop bounded operators}. 
These bounds can be established by employing Nevanlinna's theory of value distribution
similar to what we did in \cite{MyUrb08, MyUrb10.1}.
The main tool we use is Nevalinna's first main theorem (FMT) which we now describe briefly.
There are several complete accounts of it in the literature, for example in 
\cite{Nevbook74, Nevbook70, CherryYe01, GO08}.

The theory of value distribution of a meromorphic function $f:\C\to \cbar$ relies on some naturally to $f$  associated functions for which we use standard notations. For example, $n(r,w)$ or $n_f(r,w)$ is used for the 
counting function which desxribes the number of $w$--points (counted with multiplicity) of modulus at most $r$.
The average or integrated counting number $N(r,w)$ is related to $n(r,w)$ by $dN(r,w)/dr= n(r,w)/r$.

Concerning the characteristic function $\sT(r)=\sT_f (r)$ of $f$, we use the Ahlfors-Shimizu spherical version of it 
which measures the average
covering number of the Riemann sphere of the restriction of $f$ to the disk of radius $r$:
\beq \label{5.2 1}
\sT(r)= \int_0^r\left( \frac1\pi \int \int _{|x+iy|\leq t}\frac{|f'|^2}{(1+|f|^2)^2} dxdy \right)\, \frac{dt}{t}
=\int _0^r A_f(t)\frac{dt}{t} \,.
\eeq
The exponential growth of this function determines the order $\rho (f)$ of $f$ since we have
$$\rho (f) = \limsup _{r\to \infty} \frac{\sT(r)}{r} \, .$$

Nevanlinna's first main theorem (FMT) as stated in \cite{Eremenko00} (see also \cite{CherryYe01, GO08}) yields:
\bthm\label{fmt}
Let $f:\C\to \cbar$ be meromorphic of finite order. Then, with the notations above,
$$N(r,w) \leq \sT(r) + \log \frac{1}{[f(0) ,w]} \quad \text{for every }\; w\in \cbar \text{ and } r>0$$
where $[a,b]$ denotes the chordal distance on the Riemann sphere (with, in particular, $[a,b]\leq 1$, $a,b\in \cbar$).
\ethm


\

\noindent
{\bf Proof of Proposition \ref{prop bounded operators}.}
Remember first that we have the normalization Condition \ref{C2} and
thus \eqref{appendix label 1}:  $0\in \J_x$ and $|f_x(0)|\leq T$, $x\in X$.

Secondly, by \eqref{2.?} along with the Remark \ref{3.3 1}, 
$$\pfx \1(w)\leq \frac{\kappa^t  }{(1+|w|)^{(\al_2-\tau ) t}}
\sum_{f_x(z)=w}\big( 1+|z|\big)^{-t\hat\tau}\quad \text{for every} \quad 
w\in  \cV_{\d_0} ( \jul _{\shift (x)} )\;.$$
Combined with the distortion Lemma \ref{lem dist operator} it follows that the
required estimations follow if there exists $C>0$ such that
$$\sum_{f_x(z)=w}\big( 1+|z|\big)^{-t\hat\tau}\leq C \quad \text{for every} \quad 
w\in  \cV_{\d_0} ( \jul _{\shift (x)} )\setminus \D(f_x(0) , \d_0 /2)\;, \;\; x\in X \;.$$
Observe that
 $$\sum_{f_x(z)=w}\big( 1+|z|\big)^{-t\hat\tau}\leq \sum_{f_x(z)=w} \max \{1,  |z|\}^{-t\hat\tau}
 = n_{f_x} (1,w) + \sum_{f_x(z)=w\, ,\;  |z|>1} |z|^{-t\hat\tau}\,.$$
  The second term can be treated 
 by means of two integrations by part and an application of Theorem \ref{fmt} (this is completely standard, compare  also \cite[p.16]{MyUrb10.1}):

\begin{align*}
\sum _{\ba{c} f_x(z)   = w\\ \, |z|>1\ea} |z|^{-t\hat\tau} & = \int_{1}^\infty \frac{d\, n_{f_x}(r,w)}{r^{t\hat\tau}}
=-n_{f_x}(1,w) +{t\hat\tau}\int_{1}^\infty \frac{ n_{f_x}(r,w)}{r^{{t\hat\tau}+1}}dr \\
&\leq -n_{f_x}(1,w)-{t\hat\tau}N_{f_x}(1,w)
+({t\hat\tau})^2\int_{1}^\infty \frac{ N_{f_x}(r,w)}{r^{{t\hat\tau}+1}} dr \\
&\leq -n_{f_x}(1,w) +
({t\hat\tau})^2\int_{1}^\infty \frac{\sT_{f_x}(r) }{r^{{t\hat\tau}+1}} dr
+  ({t\hat\tau})^2\log \frac{1}{[f(0) ,w] }\int_{1}^\infty \frac{dr}{r^{{t\hat\tau}+1}}\,.
\end{align*}
Since $\sT_{f_x}(r) \leq C_\rho r^\rho$ (Condition \ref{C1}),
$$\sum_{f_x(z)=w}\big( 1+|z|\big)^{-t\hat\tau}\leq ({t\hat\tau})^2
\frac{C_\rho}{\hat \tau t - \rho }
+  {\hat\tau}t\log \frac{1}{[f(0) ,w] } \,.$$
The second term is uniformly bounded since we assumed $|w-f_x(0)|\geq \d_0 /2$
and since we know that $|f_x(0)|\leq T$. 
The proof is complete.

 $\hfill \square$

\subsection{Uniform second main theorem (SMT)}
Our construction of conformal measures relies on the SMT of Nevanlinna along with good estimates of the error term appearing in it. The later has been extensively studied in the 80's and 90's and the book \cite{CherryYe01} by Cherry and Ye
is an excellent reference for this topic. In particular, Chapter 2 of this book fits perfectly well to what we are doing.
The following result is a straightforward adaption of a particular case of Theorem 2.8.5 in  \cite{CherryYe01}.
We use here and throughout the whole section the notations of this book.

\bthm \label{usmt}
Let $L\geq 1$ and set $b_1=b_1(L)= e(1+ (Le^e)^2)$ and $r_0=r_0(L)=Le^e$.
Let $\rho >0$ and $C_\rho >0$. Then, for every non--constant meromorphic function $f:\C\to \cbar$
and every three distinct points $a_1,a_2,a_3\in \cbar$ verifying
\begin{itemize}
\item[(1)] \$ \frac1L \leq f^\# (0) = \frac{|f'(0)|}{1+|f(0)|^2}\leq L$,
\item[(2)] \$ f(0)\not\in \{a_1,a_2,a_3\}$ and
\item[(3)] \$ \sT_f (r) \leq C_\rho r^\rho$, $r>0$, the following holds:
\end{itemize}
$$\sum_{j=1}^3 N_f(a_j,r) \geq \sT_f(r) -S(r,a_1,a_2,a_3) \quad \text{for every} \quad r\geq r_0$$
where
\begin{align*}
S(r,a_1,a_2,a_3) & = 2\log (108 + 18\log 2)+ \frac12\log b_1+1 +4\log \sT_f(r)\\
& + \left(\frac32(\rho-1)+\frac12\right)\log r + \log L +\mathring D (a_1,a_2,a_3)\\
& \leq  b_6+ 6\rho \log r + \mathring D (a_1,a_2,a_3)\, ,
\end{align*}
 $\mathring D (a_1,a_2,a_3)=-\log \prod_{i\neq j} [a_i,a_j] +2\log 2$, $[a_i,a_j]$ being the chordal distance, and where
the constant $b_6$ does depend on $L, C_\rho $ only.
\ethm

\noindent
This, in fact uniform, version of the SMT deserves some comments.

First of all, the radius $r_0$ normally depends on the function $f$ since it is chosen in order to have 
$\sT_f(r)\geq e$. However, as it is explained in Proposition 2.8.1 of \cite{CherryYe01}, if $f$ is any meromorphic function with 
\beq\label{4.2 9}
f^\#(0) \geq \frac1L\eeq
 then $\sT_f(r)\geq \log r - \log L$. Consequently, given $L\geq 1$, there exists $r_0=r_0 (L)$
such that the above SMT does hold for every $f$ that satisfies \eqref{4.2 9}. 
Inspecting the proof of Proposition 2.8.1 of \cite{CherryYe01} gives the precise number $r_0$
indicated in the above theorem.

Various formulations of the SMT and especially the ones in Chapter 2 of \cite{CherryYe01}
 involve two functions, a \emph{Khinchin function} $\psi$
and an auxiliary function $\phi$. Their role is to optimize the error term $S(r,a_1,a_2,a_3)$ 
often by the cost of a larger \emph{exceptional set} $E$, i.e. set of radii $r\geq r_0$
such that SMT does only hold if $r\not\in E$ and this set satisfies
\beq \label{5.2 2}
\int_ E\frac{dr}{\phi (r)}\leq 2 k_0 (\psi ) = 2 \int _e^\infty \frac{dx}{x\psi (x)}\, .
\eeq
For our application we do not care about 
a minimal error term and thus  we did a more or less arbitrary chose $\psi (x)=x$. 
We equally well could have made Nevanlinna's choice $\psi (x)=(\log x)^{1+\ep}$. But our choice leads to
a nicer expression of the error term.

The choice of $\phi$ is more subtle since we need the SMT estimation for every $r\geq r_0$.
 A precise argument how to remove the exceptional set is in Nevanlinna's book \cite[p. 257]{Nevbook74}
and it is only possible since we deal with functions that have finite order. Indeed, the assumption (3)
implies that the order $\rho(f)\leq \rho$ and that the variation of the characteristic function is bounded in the following way. From the definition of $\sT_f$ in \eqref{5.2 1} follows that
$A_f(r)\leq \int _r^{er}A_f(t)\frac{dt}{t} \leq \sT_f (er) \leq C_\rho (er)^\rho  \,.$
Therefore, if $r_0\leq r_1< r_2$ then
$$\sT_f(r_2) - \sT_f(r_1) = \int _{r_1}^{r_2} A_f(t)\frac{dt}{t} \leq C_\rho e^\rho \rho^{-1} \left( r_2^\rho -r_1^\rho\right)\,.$$
Choose now, and that what we did in the above SMT, the function $\phi (r)=r^{-(\rho -1)}$.
If the interval $(r_1, r_2)\subset E$ then it results from \eqref{5.2 2} that this variation is bounded
$$\sT_f(r_2) - \sT_f(r_1)\leq C_\rho e^\rho 2k_0(\psi )$$ 
and from this it is not hard to see how to remove the exceptional set.

\vfill
\pagebreak

\bibliographystyle{plain}
\bibliography{/Users/volkmay/Math/biblio_VM}


\end{document}